# PERIODICITY IN THE TRANSIENT REGIME OF EXHAUSTIVE POLLING SYSTEMS

By I. M. MacPhee, M. V. Menshikov, S. Popov[1] and S. Volkov

*University of Durham, University of Durham, University of São Paulo and University of Bristol*

We consider an exhaustive polling system with three nodes in its transient regime under a switching rule of generalized greedy type. We show that, for the system with Poisson arrivals and service times with finite second moment, the sequence of nodes visited by the server is eventually periodic almost surely. To do this, we construct a dynamical system, the triangle process, which we show has eventually periodic trajectories for almost all sets of parameters and in this case we show that the stochastic trajectories follow the deterministic ones a.s. We also show there are infinitely many sets of parameters where the triangle process has aperiodic trajectories and in such cases trajectories of the stochastic model are aperiodic with positive probability.

**1. Introduction.** A polling system has $N$ nodes where jobs arrive and queue and a single server which switches between the nodes to process the jobs. In exhaustive polling systems the server processes all jobs at its current node $i$, say, including any that arrive while jobs there are being processed, before switching to another node $j$, chosen by some rule. Conditions for transience/recurrence of polling systems using a *greedy* switching rule were given in [5, 6]. There is a critical case which has been investigated in [8] and [9]. In this paper we show that, for an exhaustive polling system with $N = 3$ nodes, with arrival streams and service times putting the system in its transient regime and switching according to a *threshold rule* (a type of generalized greedy rule), the sequence in which the server visits the nodes is eventually periodic for almost all choices of threshold parameters.

---

Received May 2004; revised November 2005.
[1]Supported in part by CNPq (302981/02–0) and FAPESP (04/03056–8).

To show this, we consider in Section 2 the embedded Markov chain in $\mathbf{Z}_+^N$ where we observe the polling system at service and switching time completions. We construct a dynamical system using the vector field of expected drifts of the chain for the server at each node. As the polling system is transient, this dynamical system exits any finite ball eventually, so we project it onto the unit simplex. With $N = 3$ nodes, we call this the *triangle process*, as it lives in a triangle—it has piecewise linear trajectories that change direction when they meet the triangle boundary so it looks something like a billiards model, but where "reflections" are caused by changes of dynamics due to the server switching to another node. We call any trajectory which returns to its start point after finitely many server switchings an *orbit*.

It is worth mentioning that this deterministic system closely resembles (and in some cases *is*) the so-called affine interval exchange transformation model. However, the affine interval exchange transformation here is contracting so it is not a bijection and it reverses orientation so it is not order-preserving. It is thus rather different from the usual interval exchange transformation model that has been much studied during the past 30 years. We comment more on this in Section 5.

We state our main results in Section 3. Theorem 3.1 states that for a.s. all choices of switching thresholds the triangle process has a finite number of orbits and every trajectory converges toward one of these orbits. We say the triangle process is stable in this case. That there are no more than four orbits for any triangle process parameters is the content of Theorem 3.2. In Theorem 3.3 we show that there are (infinitely many) choices of the switching thresholds that lead to the existence of aperiodic trajectories of the triangle process. These results are proved in Sections 4 and 5. Using these results for the deterministic triangle process, we show, in Theorem 3.4, that the projections of the trajectories of the stochastic queueing model a.s. converge onto the orbits of the triangle process, when it is stable. This implies the periodicity of the sequence of nodes receiving service for the polling system. Finally, in Theorem 3.5, we show that when the triangle process has nonperiodic trajectories then, in some situations, the sequence of nodes visited by the stochastic model is periodic with positive probability and aperiodic with positive probability. These results are proved in Sections 6 and 7.

**2. System description.** An exhaustive polling system has $N$ nodes where jobs queue and a single server which switches to the next node $j$, chosen with some rule, after processing all jobs at the current node $i$, including any that arrive while jobs at $i$ are being processed. In the general model switching takes a time which depends upon the pair $i$, $j$. We will assume the following at each node $i$: the arrival processes are independent Poisson streams with arrival rate $\lambda_i$; the service times are i.i.d. and independent of the arrivals with finite mean $\mu_i^{-1}$ and variance $\sigma_i^2 < \infty$; the switching times are independent with finite means that depend upon the initial and final node while the server just waits for the next arrival at any node when it completes service at $i$ and finds the system is empty.

For the results in our paper, we will assume that the switching times are zero. We believe our argument can be extended to a more general situation, but our main interest is in the transient case, where the system behavior is not sensitive to the switching time distributions, though of course it is affected by the switching rule.

The methods we know for showing that periodicity of the node sequence for the deterministic triangle process implies the same for the stochastic process are applicable only to Markov processes. As the sequence of nodes visited is determined by the jump chain, we will consider the *discrete* time process

$$\Xi = \{(\xi(t); s(t))\}, \qquad t = 0, 1, \ldots,$$

where $\xi(t)$ describes the queue lengths and $s(t)$ the server location at the epochs of service time and switching time completions. Its state space is $\mathbb{Z}_+^N \times \{1, 2, \ldots, N\}$. That $\Xi$ is well defined, irreducible and aperiodic follows from the assumptions of Poisson arrivals and finite first moments of the service times. For later convergence arguments, we also require finite second moments of the service times and under these conditions, standard results imply this embedded chain is essentially equivalent to the continuous time process as regards transience/recurrence. The transition probabilities can be computed via the Laplace transforms of the service times, but we do not need them explicitly at any point. A standard conditioning argument readily produces the expected one-step mean drifts

$$\mathbf{E}(\xi_i(t+1) - \xi_i(t) | (\xi(t); s(t)) = (x; j)) = \lambda_i \mu_j^{-1} - I_{\{i=j\}}, \qquad i = 1, \ldots, N, \tag{1}$$

for $x \in \mathbb{Z}_+^N$ with $x_j \geq 1$ and the server at node $j$ and we make considerable use of these. A similar result holds for expected drifts during switching times, but we will not detail these as, in fact, we will assume switching to occur instantaneously. Other arrival and service processes for which there is an embedded Markov chain can be found, but their treatment needs no significant extension of our methods.

We will start by stating the known explicit conditions for recurrence/transience. For exhaustive polling models, the conditions for recurrence/transience depend only upon the *total loading* (or traffic intensity)

$$\rho = \sum_{i=1}^{N} \rho_i \qquad \text{where } \rho_i = \lambda_i/\mu_i \text{ at each node } i$$

for a great many (and seemingly all *sensible*) switching rules under the assumption that switching times have finite first moment. We have the following:

THEOREM 2.1. *The process $\Xi$ is positive recurrent if $\rho < 1$, transient if $\rho > 1$.*

REMARKS. Foss and Last [5] establish the result for a more general model. The method of proof is via Lyapunov functions as described in [1, 4] or [10]. The papers by MacPhee and Menshikov [8] and Menshikov and Zuyev [9] consider the critical case $\rho = 1$, where the behavior of the system depends strongly on the first two moments of the switching time distribution. Foss and Last [6] obtain the result of this theorem for nonexhaustive polling systems under a greedy switching policy.

*Generalized greedy switching rules.* When the server has completed all the tasks at its current node, it chooses its next node using a switching rule which we will assume depends upon the queue lengths, that is, the rule is a function $\mathfrak{R} : \partial\mathbb{R}_+^N \to \{1, 2, \dots, N\}$, where $\partial\mathbb{R}_+^N = \bigcup_j \{y \in \mathbb{R}_+^N : y_j = 0\}$. Any $\mathfrak{R}$ must satisfy $\mathfrak{R}(y) \neq j$ if $y_j = 0$ and a variety of such rules have been studied in the literature. We will study only a class of generalized greedy rules defined as follows. For each node $j$, there is a vector of positive weights $b_j = (b_{j1}, \dots, b_{jN})$ and at states $y \in \partial\mathbb{R}_+^N$ with $y_j = 0$, $\mathfrak{R}(y) = i$, where $b_{ji} y_i = \max_k (b_{jk} y_k)$ (for our results it is not important how ties are resolved). The simple greedy rule is the special case where all $b_{jk}$ are equal.

*A deterministic model.* Our subsequent analysis of the transient case is based around the following deterministic model of the system. Consider a particle moving in $\mathbb{R}_+^N \times \{1, \dots, N\}$ with linear dynamics given by the one-step mean drifts of $\Xi$ as calculated in (1). With its position denoted $y = (y_1, y_2, \dots, y_N)$ and the server at node $j$, we see that at $(y, j)$ with $y_j > 0$ the particle has velocity

$$\mu_j^{-1} \sum_{i=1}^{N} \lambda_i e_i - e_j = \mu_j^{-1} \left( \sum_{i \neq j} \lambda_i e_i + (\lambda_j - \mu_j) e_j \right),$$

where the $e_i$ denote the axial unit vectors in $\mathbb{R}^N$. If a point with $y_j = 0$ is reached, then a different set of dynamics (corresponding to the server switching to another node) is chosen instantaneously according to some generalized greedy switching rule $\mathfrak{R}$.

The trajectory $y(t)$ of our particle is constructed as follows. From start point $y(t_0) = \bar{y} \in \mathbb{R}_+^N$ with $\bar{y}_j > 0$ and the server at node $j$, the particle

travels along the line

$$\ell_j(\bar{y}) = \left\{ y \in \mathbb{R}^N : y = \bar{y} + \mu_j^{-1}(t - t_0)\left(\sum_{i \neq j} \lambda_i e_i + (\lambda_j - \mu_j)e_j\right), t \in \mathbb{R} \right\}$$

through $\bar{y}$ with the appropriate velocity. If $\rho_j \geq 1$, the particle never reaches $\partial\mathbb{R}^N_+$, while if $\rho_j < 1$, the particle reaches

$$y(t_1) = \sum_{i \neq j} \left( \bar{y}_i + \frac{\lambda_i \bar{y}_j}{\mu_j - \lambda_j} \right) e_i \in \partial\mathbb{R}^N_+ \tag{2}$$

at time $t = t_0 + \bar{y}_j/(1 - \rho_j) =: t_1$. Now the server switches to node $\mathfrak{R}(y(t_1))$ and the next and subsequent pieces of the trajectory are computed as above.

*The node process.* In the transient case the trajectory exits from any finite ball eventually, so, in order to study whether there is any periodicity in the order the server visits the queues, we project onto the unit simplex. Project the lines $\ell_j(\bar{y})$ onto the hyperplane $\mathcal{S}_1 = \{y \in \mathbb{R}^N : \sum_i y_i = 1\}$ using

$$\Lambda : \mathbb{R}^N \setminus \{0\} \to \mathcal{S}_1 \qquad \text{where } \Lambda(y) = \frac{y}{\sum_i y_i},$$

which maps each line $\ell_j(\bar{y})$ onto the intersection of $\mathcal{S}_1$ with the plane containing $\ell_j(\bar{y})$ and the origin. The feasible positions for the particle will be mapped onto points in the unit simplex $\mathcal{S}_1^+ = \{y \in \mathbb{R}^N_+ : \sum_i y_i = 1\}$. For starting points $\bar{y}$ and $\alpha\bar{y}$ for any $\alpha > 0$, the lines $\ell_j(\bar{y})$ and $\ell_j(\alpha\bar{y})$ have exactly the same image in $\mathcal{S}_1$, so we can restrict our attention to reference points $\bar{y} \in \partial\mathcal{S}_1^+$, the boundary of $\mathcal{S}_1^+$.

Under the condition $\rho_j > 1$, the server in the stochastic process can remain serving at queue $j$ indefinitely, so we will only consider the cases where $\rho_j < 1$ at all queues $j$. Under the condition $\rho_j < 1$, we see from equation (2) that the trajectory $y(t)$ leaving $\bar{y} \in \partial\mathcal{S}_1^+$ along $\ell_j(\bar{y})$ next reaches $\partial\mathbb{R}^N_+$ in finite time at the point

$$\sum_{i \neq j} \left( \bar{y}_i + \frac{\lambda_i \bar{y}_j}{\mu_j - \lambda_j} \right) e_i \qquad \text{with projection } \sum_{i \neq j} \frac{(\mu_j - \lambda_j)\bar{y}_i + \lambda_i \bar{y}_j}{(\mu_j - \lambda_j) + \mu_j \theta_j \bar{y}_j} e_i \in A_j^0,$$

where $A_j^0 = \{y \in \mathbb{R}^N : \sum_i y_i = 1, y_j = 0, y_i \geq 0 \text{ for } i \neq j\}$ for $j = 1, 2, \ldots, N$ are regions of the switching boundary $\partial\mathcal{S}_1^+$ and

$$\theta_j := \mu_j^{-1} \sum_i \lambda_i - 1, \qquad j = 1, 2, \ldots, N. \tag{3}$$

There is a very nice geometric description of the projected process. For any fixed $j$, the lines $\ell_j(\bar{y})$ through different points $\bar{y}$ are parallel so they are concurrent after projection, that is, their image lines on $\mathcal{S}_1$ are either parallel or share a common focus, $v_j$, say. The line $\ell_j(\bar{y})$ with $\bar{y} = \sum_{i \neq j} \lambda_i e_i + (\lambda_j - \mu_j) e_j$ passes through 0 so $v_j$ is the point where $\ell_j(\bar{y})$ meets $\mathcal{S}_1$ if this happens. This is the case whenever $\theta_j \neq 0$, in which case

$$v_j = \frac{1}{\mu_j \theta_j} \left( \sum_{i \neq j} \lambda_i e_i + (\lambda_j - \mu_j) e_j \right). \tag{4}$$

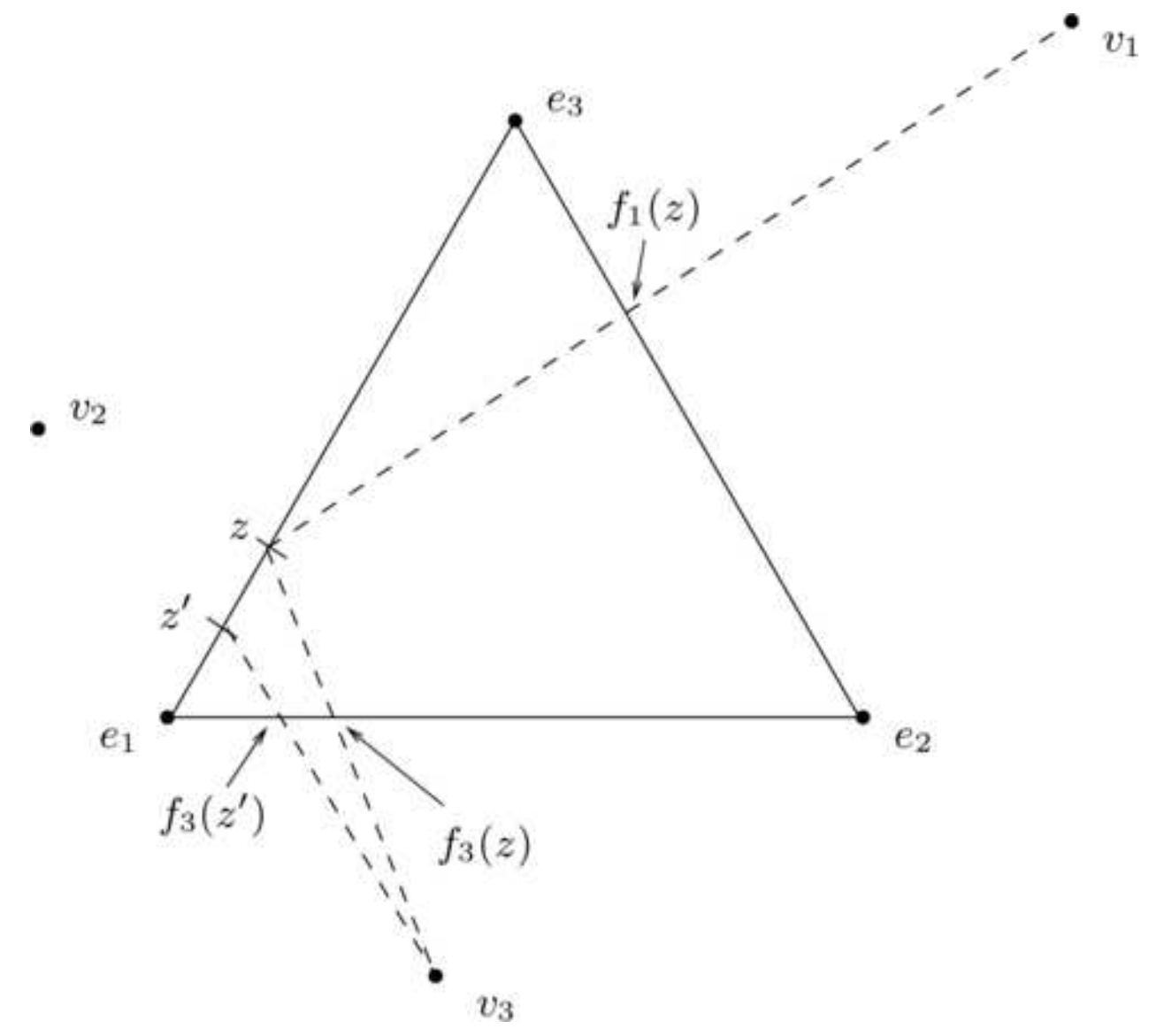

FIG. 1. *The node process when $N = 3$.*

When $\theta_j = 0$, the $\ell_j(\bar{y})$ with $\bar{y} \in \mathcal{S}_1$ all lie in $\mathcal{S}_1$ and the projected lines through different $\bar{y}$ are parallel.

The focus point $v_j$ will be outside $\mathcal{S}_1^+$ when $\theta_j < 0$ or when $\theta_j > 0$ and $\rho_j = \lambda_j / \mu_j < 1$, conditions which correspond to the next switching event taking place after some finite time.

As our aim is to study the switching sequence, we assume from now on that

$$\rho_i < 1, \qquad i = 1, \dots, N \quad \text{and} \quad \rho = \sum_{i=1}^{N} \rho_i > 1 \tag{5}$$

(which means that the whole system is transient, yet the server does not get stuck at any individual node). We can now define the discrete dynamical system $Z = \{z(t)\}, t = 0, 1, \dots,$ living on $A^0 \equiv \bigcup_i A_i^0 = \partial \mathcal{S}_1^+$. Define mappings

$$f_j(z) = \sum_{i \neq j} \frac{(\mu_j - \lambda_j) z_i + \lambda_i z_j}{(\mu_j - \lambda_j) + \mu_j \theta_j z_j} e_i \in A_j^0, \qquad z \in A^0 \setminus A_j^0, \ j = 1, \dots, N. \tag{6}$$

The image under $f_j$ of a point $z \in A_i^0$ lies at the intersection of the line through $z$ and $v_j$ and the boundary region $A_j^0$, as shown for the case $N = 3$ in Figure 1. For given $z(0) \in \partial\mathcal{S}_1^+$ let

$$z(t+1) = \varphi(z(t)), \qquad t = 0, 1, \dots, \tag{7}$$
$$\text{where } \varphi(z) = \sum_{j=1}^{N} I_{\{\Re(z)=j\}} f_j(z), z \in A^0.$$

We will call $Z$ the *node process*, as it records information about the projection of the dynamical system $y(t)$ only at switching epochs.

**3. The three node case.** From this point our analysis is restricted to the case where the system has three nodes and the switching decision is made using a generalized greedy rule. This enables us to give clear statements of our results and methods, but is complex enough to be very interesting in our opinion.

Recall that when the server has completed all jobs at queue $i$, the switching rule is defined by parameters $b_{ij} > 0$ and selects the next node $k$ when $b_{ij}x_j < b_{ik}x_k$ for $i, j, k$ any permutation of 1, 2, 3. Projection by $\Lambda$ onto $\mathcal{S}_+^1$ reduces the switching boundaries $\{x : b_{ij}x_j = b_{ik}x_k,\ x_i = 0\}$ to *decision points*, one on each side of the triangle $A^0$ as shown in Figure 2. This rule can be applied to the stochastic process $\Xi$ and the dynamical process $y(t)$ as stated with the following resolution of boundary cases. For the stochastic process $\Xi$, a randomized rule may be used whenever $b_{ij}x_j = b_{ik}x_k$, so we will consider the consequences for the node process $Z$ of both possible decisions at such points. Specifically, we will consider trajectories of the node process which branch when they exactly hit the decision points. From now on, we will call the node process $Z$ the *triangle process* and specialize our notation. The construction in Section 2 is somewhat abstract, but the system which we have defined and wish to study is really very simple to describe. Figure 2 shows a trajectory starting from $z(0) = z$ near $e_3$. This is mapped to $z(1)$ along the line from $z$ to $v_3$ and then to $z(2)$ just below $d_1$. As $z(2)$ is toward $e_2$ from $d_1$, the next point $z(3)$ is on the line from $z(2)$ to $v_2$. From here the trajectory continues toward $v_3$, then $v_1$, then back toward $v_3$ and from there toward $v_2$ and so on.

We note that the mapping which determines the trajectories is not continuous at the $d_i$. However, during extensive numerical investigation, we found that, for all the configurations of rates and decision points we tried, the trajectories we examined converged toward periodic orbits.

*Triangle process notation.* Let $\hat{\imath}, \hat{\jmath}, \hat{k}$ denote the values 1, 2, 3 or either of its *cyclic* permutations 2, 3, 1 or 3, 1, 2 and let $i, j, k$ denote any permutation of 1, 2, 3. For $z=(1-x)e_{\hat{\jmath}}+xe_{\hat{k}}\in A^0_{\hat{\imath}}$ with $x\in[0,1]$, we define $\mathrm{Side}(z):=\hat{\imath}$, $\pi(z):=x$ and write $z=(x,\hat{\imath})$ (with standard right-hand axes, increasing $x$ corresponds to going around the triangle anticlockwise).

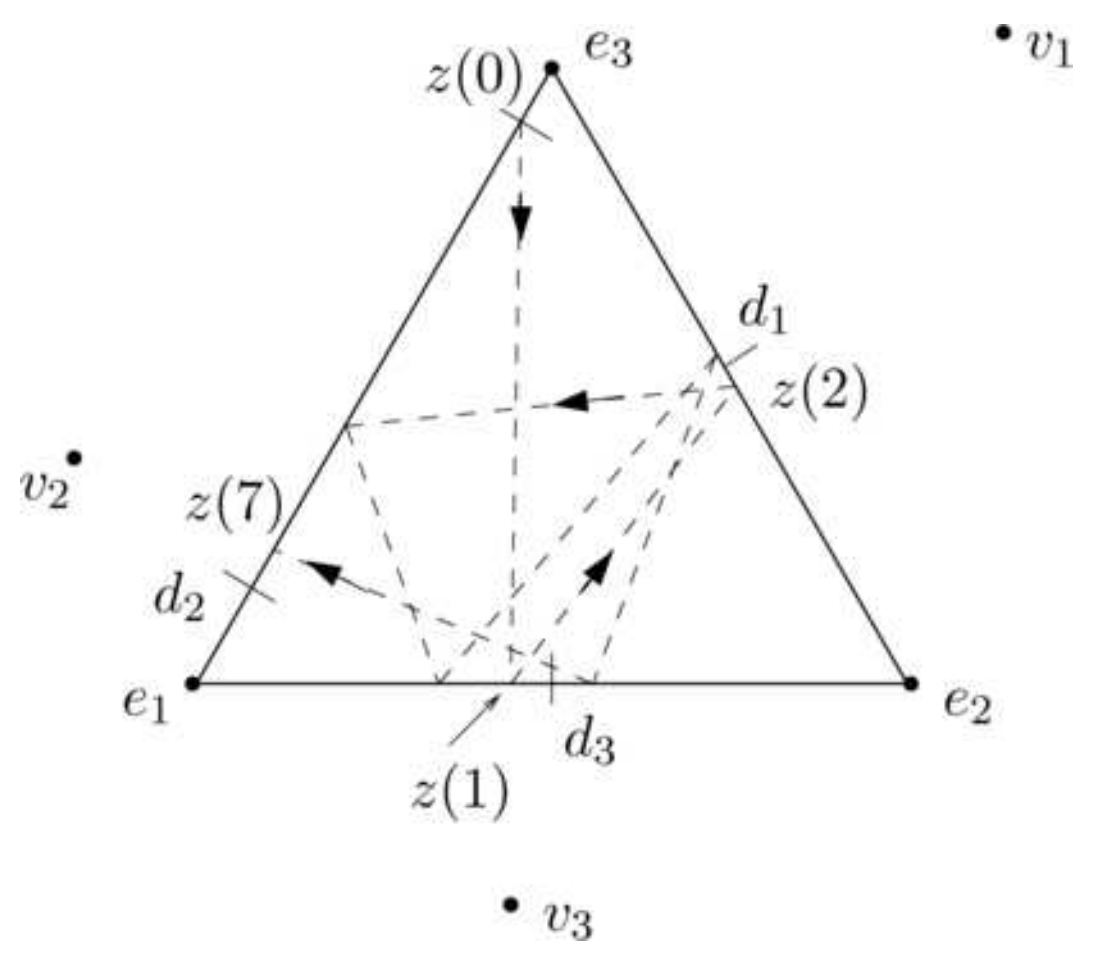

FIG. 2. *The triangle process trajectory.*

For the triangle process $Z$, any generalized greedy rule $\mathfrak{R}$ has the form

$$\mathfrak{R}(x,\hat{\imath})=\begin{cases}\hat{\jmath}, & x\le d_{\hat{\imath}},\\ \hat{k}, & x\ge d_{\hat{\imath}},\end{cases}$$

for values $d_{\hat{\imath}}=(1+b_{\hat{\imath}\hat{k}}/b_{\hat{\imath}\hat{\jmath}})^{-1}\in(0,1)$ for $\hat{\imath}=1,2,3$. The *decision points* $(d_i,i)$ will usually be written simply $d_i$. We will usually refer to rules of this type as *threshold* rules.

For $z=(x,i),z'=(x',j)$, we will use distance $|z-z'|_1:=\|z-z'\|/\sqrt{2}$, the Euclidean distance scaled so that when $i=j,|z-z'|_1=|\pi(z)-\pi(z')|=|x-x'|$, that is, it is length along the side of the triangle.

The forward mapping $\varphi$ defined in (7) is 1–1 except at the $d_i$ which have two images, while the inverse of $\varphi$ is 1–1 where it exists. Call $z$ a *pre-image* of $z'$ when $z'=\varphi^{(t)}(z)$ for some $t\ge 1$.

DEFINITION. A sequence $z(t),t=0,1,\dots,$ satisfying (7) is a *trajectory* of $Z$. If $z(0)$ is a pre-image of a decision point $d_i$, then there are at least two trajectories starting from $z(0)$.

We will study closely the sets $A^t=\varphi(A^{t-1}),t=1,2,\dots,$ which are such that $A^t$ contains the possible locations for $z(t)$ from all initial points $z(0)\in A^0$. Alternatively, $A^0\setminus A^t$ is the set of points with fewer than $t$ pre-images, a description which we will use later. The set $A^1=\varphi(A^0)$ is a strict subset of $A^0$ for any decision points since, for example, $e_i\notin A^1,i=1,2,3$. Hence, $A^0 \supset A^1 \supset \cdots \supset A^t$ for $t = 2, 3, \ldots.$ Further, $A^1$ is a union of three disjoint closed intervals, one in each $A_i^0$. We introduce now the notation $A_i^t \equiv A^t \cap A_i^0$ for later use—$A_1^1$ is depicted in Figure 3. As only intervals containing a decision point will be split by $\varphi$, it follows that $A^t$ is a union of at most $3t$ disjoint closed intervals.

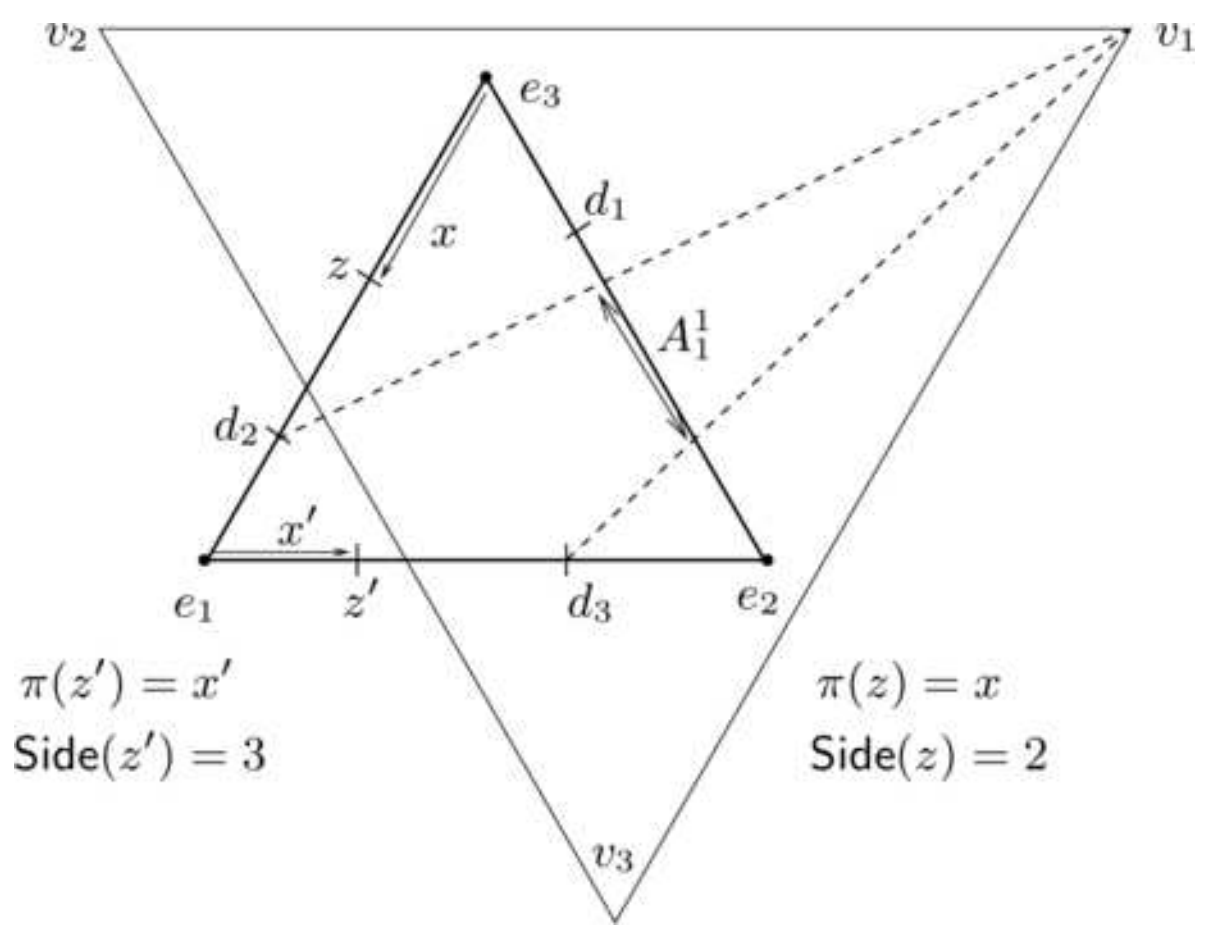

FIG. 3. *The triangle process.*

*Periodicity definitions.* We call a trajectory $z(t)$, $t = 0, 1, \ldots,$ *eventually-$m$-periodic* if there is $N > 0$ such that, for all $n \geq N$,

$$\mathrm{Side}(z(n+m)) = \mathrm{Side}(z(n)).$$

A trajectory which is not eventually-$m$-periodic for any $m$ is called *nonperiodic.*

We say that $Z$ has a *periodic orbit* if there is a finite sequence of points $u_1, u_2, \ldots, u_m$ with $\phi(u_i) = u_{i+1}$ for $i = 1, 2, \ldots, m-1$ and $\phi(u_m) = u_1$. The sequence $\mathrm{Side}(u_1), \mathrm{Side}(u_2), \ldots, \mathrm{Side}(u_m)$ is the *node-cycle* of the orbit. We will consider cyclic permutations of orbits/node-cycles to be equivalent to the original orbit/node-cycle. All orbits are disjoint since, for any point $u$ on an orbit, $\varphi^{(-t)}(u)$ is uniquely defined for all $t \geq 0$ so there is no point at which two orbits could join.

A trajectory *converges onto an orbit* $u_1, \ldots, u_m$ when there exists $n_0 \geq 0$ such that

$$|z(mt + n + n_0) - u_n|_1 \to 0 \qquad \text{as } t \to \infty \text{ for each } n = 1, 2, \ldots, m,$$

that is, for each $n$, the subsequence of the trajectory corresponding to phase $n$ of the node-cycle converges to $u_n$. For any point $z = (x, i) \in A_i^0$ and $0 < \varepsilon < x$, let

$$\mathcal{N}_\varepsilon(z) = \{z' \in A_i^0 : |\pi(z') - x|_1 < \varepsilon\} \tag{8}$$

denote the $\varepsilon$-neighborhood of $u$, $\mathcal{N}_\varepsilon^+(z) = \{z' \in A_i^0 : 0 \le \pi(z') - x < \varepsilon\}$ and $\mathcal{N}_\varepsilon^-(z) = \{z' \in A_i^0 : 0 \le x - \pi(z') < \varepsilon\}$ the one-sided $\varepsilon$-neighborhoods of $z$. We will say an orbit is *stable* when, for each $u_n$, there exists $\varepsilon > 0$ such that trajectories starting from any $z \in \mathcal{N}_\varepsilon(u_n)$ converge onto the orbit. When trajectories started from $z \in \mathcal{N}_\varepsilon^+(u_n)$ but not from $\mathcal{N}_\varepsilon^-(u_n)$ (or vice versa) converge onto the orbit, we say the orbit is *stable on one side* and otherwise we say the orbit is *unstable*. An orbit containing a decision point must be unstable if it is of odd length $m$ and may be stable on one side if $m$ is even.

We are now in position to state our main results for the triangle process $Z$. These hold for any set of parameters $\lambda_i$, $\mu_i$, $i = 1, 2, 3$, satisfying conditions (5).

THEOREM 3.1. *For almost all decision points $d_i \in A^0, i = 1, 2, 3$, the triangle process $Z$ has finitely many periodic orbits. For such sets of decision points, all trajectories $z(t)$ are eventually periodic and each converges onto one of these orbits as $t \to \infty$.*

The quantifier "almost all" is used in the sense of Lebesgue measure $\times$ counting measure on $[0, 1] \times \{1, 2, 3\}$. Due to the convergence behavior of all trajectories, we will call this the *stable case*. Typically, in this case all the orbits are stable, but unstable orbits are possible when an orbit includes a decision point.

THEOREM 3.2. *For any set of decision points $d_i, i = 1, 2, 3$, there are at most four periodic orbits.*

The next result shows that not all choices of decision points result in all trajectories being eventually periodic.

THEOREM 3.3. *There is an uncountable set of decision points for which $Z$ has nonperiodic trajectories.*

These results have important consequences for the behavior of the underlying stochastic process $\Xi = \{\xi(t)\}_t$. We now introduce the random times $\tau_t, t = 1, 2, \ldots$, at which the server changes queues. Under the assumptions in (5) on the parameters, the $\tau_t$ are all a.s. finite. The next result concerns the stochastic triangle process

$$\zeta(t) = \Lambda(\xi(\tau_t)), \qquad t = 1, 2, \ldots. \tag{9}$$

This process has trajectories living on $A^0$, so we can use the same definition of convergence for it as for the triangle process.

THEOREM 3.4. *Suppose the service times have variances $\sigma_i^2 < \infty$. For any set of decision points $(d_1, d_2, d_3)$ such that $Z$ is stable, the stochastic process $\zeta$ is also stable in the sense that a.s. each trajectory of $\zeta$ converges onto one of the periodic orbits of $Z$.*

THEOREM 3.5 (No zero–one law). *There exist configurations of the decision points such that the stochastic process $\zeta$ has nonperiodic trajectories with positive probability.*

We construct an example with trajectories that can converge to a periodic orbit of $Z$ with positive probability and also can be nonperiodic with positive probability.

*Some open problems.*

- It seems from numerical computations that there are, in fact, at most three orbits. Is this correct or are there examples with four orbits?
- Is it also true that, for any given set of decision points $d_i$, the sets of parameters $\lambda, \mu$ where there are nonperiodic orbits has measure zero?
- Which of the results proved here also hold when there are four or more nodes?
- Can the stochastic process $\zeta$ converge with positive probability to an orbit of even length which contains a decision point and is only stable on one side?

**4. Proofs for the triangle process.** The conditions for transience or recurrence of the stochastic process $\Xi$ are in terms of the traffic intensities $\rho_i$ in accord with the intuition that considering weighted work at nodes rather than just the numbers of queued jobs should not affect such properties. Such a re-weighting also helps simplify the treatment of the triangle process. To describe its effect, we introduce the fractional linear functions $F_\alpha : [0,1] \to [0,1]$, where, for any $\alpha > 0$, $F_\alpha(x) = x/(\alpha + (1-\alpha)x)$. These have a key composition property $F_\alpha \circ F_\beta = F_{\alpha\beta}$.

LEMMA 4.1 (Re-weighting). *Consider the triangle process $Z$ with parameters $\lambda_i$, $\mu_i$, $d_i$ for $i = 1,2,3$ and define $T : A^0 \to A^0$ by $T(z) = F_{\mu_{\hat{k}}/\mu_{\hat{\jmath}}}(z), z \in A^0_{\hat{\imath}}, \hat{\imath} = 1,2,3$. $Z$ is isomorphic to the triangle process $Z'$ with parameters $\mu'_i = 1$, $\lambda'_i = \lambda_i/\mu_i = \rho_i$ and $d'_i = T(d_i)$, $i = 1,2,3$. Specifically, for any given $z(0) \in A^0$, if $Z'$ is started from $T(z(0))$, then $z'(t) = T(z(t))$ for $t = 1,2,\ldots.$*

PROOF. We start by assuming that only parameters $\lambda_{\hat{\imath}}$, $\mu_{\hat{\imath}}$ for some $\hat{\imath}$ are transformed by multiplication by $\alpha > 0$ and briefly describe an isomorphism between trajectories of the dynamical processes $y$ and $y'$. Define $\Gamma_{\hat\imath}:\mathbb{R}_+^3 \to \mathbb{R}_+^3$ by $\Gamma_{\hat\imath}(y) = y + (\alpha-1)y_{\hat\imath}e_{\hat\imath}$, so $\Gamma_{\hat\imath}$ rescales $\mathbb{R}_+^3$ by $\alpha$ in the $e_{\hat\imath}$ direction. Let $\ell_j(\bar y)$ and $\ell_j'(\bar y)$ denote the lines through $\bar y$ parallel to the trajectories of $y(t)$ and $y'(t)$ respectively when the server is at node $j$. Then $\Gamma_{\hat\imath}(\ell_j(\bar y)) = \ell_j'(\Gamma_{\hat\imath}(\bar y))$. This space rescaling provides an isomorphism between entire trajectories whenever the switching decisions are identical, or, equivalently, $r'(\Gamma_{\hat\imath}(\bar y)) = r(\bar y)$ for all $\bar y \in \partial\mathbb{R}_+^3$.

The class of threshold policies is closed under space rescalings $\Gamma_{\hat\imath}$ with the parameters $b_{\hat\jmath\hat\imath}$ being mapped to $\alpha b_{\hat\jmath\hat\imath}$ and $b_{\hat k\hat\imath}$ to $\alpha b_{\hat k\hat\imath}$ with consequent changes to $d_{\hat\jmath}$ and $d_{\hat k}$. Now for any $\bar y \in \partial\mathbb{R}_+^3$ and $z = \Lambda(\bar y) = (x,i)$, we find that $z' = \Lambda(\Gamma_{\hat\imath}(z))$ satisfies

$$z' = \begin{cases} (x,i), & i = \hat\imath, \\ (F_{1/\alpha}(x), i), & i = \hat\jmath, \\ (F_\alpha(x), i), & i = \hat k, \end{cases}$$

where $(\hat\imath,\hat\jmath,\hat k)$ is a cyclic permutation of $(1,2,3)$. Applying successive rescalings by $1/\mu_i$ in the $e_i$ direction for $i = 1$, 2, 3 and using the composition property of $F_\alpha$ gives the result. Note that $F_\alpha$ satisfies $F_\alpha(0) = 0$, $F_\alpha(1) = 1$ and $F_\alpha'(x) > 0$, $x \in (0,1)$, so the mapping $T$ fixes the corners $e_i$ of $A^0$ and smoothly rescales the sides. $\square$

The implication of Lemma 4.1 is that we need only study the triangle process with parameters $\mu_1 = \mu_2 = \mu_3 = 1$ and general $\lambda_i$ and $d_i$ to understand the possible behavior of the general case. We will work with the rates $\mu_k = 1$ from this point. We now find, from (3), that $\theta_k = \sum_n \rho_n - 1 =: \theta > 0$ for each $k$ as we are in the transient case. It follows from standard projective geometry results that the foci $v_k$ are the vertices of an equilateral triangle $V$ with sides parallel to the $A_i^0$, as shown in Figure 3. It can also be seen directly as, from equation (4) for each $k$,

$$v_k = \frac{1}{\theta}(\rho_i e_i + \rho_j e_j + (\rho_k - 1)e_k) \quad \text{and, hence,} \quad v_i - v_j = \frac{1}{\theta}(e_j - e_i).$$

Under conditions (5), the $v_k$ are outside $\mathcal{S}_1^+$, as remarked after equation (4). Let ${}^V\!A$ denote those points of $A^0$ which are internal to $V$ and ${}^V\!A_i = A_i^0 \cap {}^V\!A$. In Figure 3, ${}^V\!A$ contains all of $A^0$ except the neighborhood of $e_1$.

When $\rho_k < \theta$ or, equivalently, $\rho_i + \rho_j > 1$, the line from $v_i$ to $v_j$ meets $A_i^0$ at the point $(1-\alpha)e_j + \alpha e_k$ with $\alpha = \rho_k/\theta \in (0,1)$, while if $\rho_i + \rho_j < 1$, then $e_k \in {}^V\!A$. If $e_{\hat k} \notin {}^V\!A$, let $J_{\hat k} = \{z \in A_{\hat\imath}^0 : x > \rho_{\hat k}/\theta\} \cup \{z \in A_{\hat\jmath}^0 : x < 1 - \rho_{\hat k}/\theta\}$ denote the corner of $A^0$ containing $e_{\hat k}$ and lying outside ${}^V\!A$. Let $J_{\hat k} = \varnothing$ when $e_{\hat k} \in {}^V\!A$ (so $J_2 = J_3 = \varnothing$, but $J_1 \neq \varnothing$ in Figure 3).

The next two lemmas state key properties of the node process mappings but ignore the effect of the decision points. Translating definition (6) into the triangle process notation, we find that $f_{\hat{\jmath}}\colon A^0_{\hat{\imath}} \cup A^0_{\hat{k}} \to A^0_{\hat{\jmath}}$ satisfies

$$(10)\qquad \begin{aligned} f_{\hat{\jmath}}(x,\hat{\imath}) &= \left(\frac{\rho_{\hat{\imath}}(1-x)}{1-\rho_{\hat{\jmath}}+\theta(1-x)}, \hat{\jmath}\right) \quad \text{and} \\ f_{\hat{\jmath}}(x,\hat{k}) &= \left(\frac{\rho_{\hat{\imath}}x+(1-\rho_{\hat{\jmath}})(1-x)}{1-\rho_{\hat{\jmath}}+\theta x}, \hat{\jmath}\right) \end{aligned}$$

and $f_{\hat{\jmath}}(e_{\hat{\jmath}}) = f_{\hat{\jmath}}(0,\hat{\imath}) = f_{\hat{\jmath}}(1,\hat{k}) = \rho_{\hat{\imath}}/(\rho_{\hat{\imath}}+\rho_{\hat{k}})$, which demonstrates continuity of the $f_j$ at the corners $e_j$ of $A^0$.

LEMMA 4.2 (Contraction property). *The function $f_{\hat{\jmath}}$ is monotone on $A^0_{\hat{\imath}} \cup A^0_{\hat{k}}$ and there exists $\gamma \in (0,1)$ and $C_{\hat{\jmath}}(\gamma) \subset A^0_{\hat{\imath}} \cup A^0_{\hat{k}}$ such that for $z$, $z' \in C_{\hat{\jmath}}(\gamma) \cap A^0_n$, $n = \hat{\imath}, \hat{k}$, we have $|f_{\hat{\jmath}}(z) - f_{\hat{\jmath}}(z')|_1 < \gamma|z-z'|_1$ and further,*

$$C_{\hat{\jmath}}(\gamma) \supset \{(x,\hat{\imath})\colon 0 \le x \le \min(1, \rho_{\hat{k}}/\theta)\} \cup \{(x,\hat{k})\colon \max(0, 1-(\rho_{\hat{\imath}}/\theta)) \le x \le 1\},$$

*so $C_{\hat{\jmath}}(\gamma) \supset {}^VA_{\hat{\imath}} \cup {}^VA_{\hat{k}}$. The containment is strict if $J_{\hat{\imath}} \cup J_{\hat{k}} \ne \varnothing$.*

PROOF. On $A^0_{\hat{\imath}}$, let $g(x) = \pi(f_{\hat{\jmath}}(x,\hat{\imath})) = \rho_{\hat{\imath}}(1-x)/(1-\rho_{\hat{\jmath}}+\theta(1-x))$. As $g'(x) = -\rho_{\hat{\imath}}(1-\rho_{\hat{\jmath}})/(\rho_{\hat{\imath}}+\rho_{\hat{k}}-\theta x)^2 < 0$, $g$ is monotone decreasing in $x$. Choose $\gamma_{\hat{\imath}}$ such that if $\rho_{\hat{k}} > \theta$, then $|g'(1)| = \rho_{\hat{\imath}}/(1-\rho_{\hat{\jmath}}) \le \gamma_{\hat{\imath}} < 1$, while if $\rho_{\hat{k}} < \theta$, $|g'(\rho_{\hat{k}}/\theta)| = (1-\rho_{\hat{\jmath}})/\rho_{\hat{\imath}} \le \gamma_{\hat{\imath}} < 1$. The same argument applies for selecting $\gamma_{\hat{k}}$ on $A^0_{\hat{k}}$ using $g(x) = (\rho_{\hat{\imath}}x + (1-\rho_{\hat{\jmath}})(1-x))/(1-\rho_{\hat{\jmath}}+\theta x)$. Finally, choose $\gamma \ge \max(\gamma_{\hat{\imath}}, \gamma_{\hat{k}})$ and set $C_{\hat{\jmath}}(\gamma) = \{z \in A^0_{\hat{\imath}} \cup A^0_{\hat{k}}\colon |g'(x)| \le \gamma\}$. □

The mapping $f_{\hat{\jmath}}$ may not be contracting near $e_{\hat{\imath}}$ or $e_{\hat{k}}$ if they are not in ${}^VA$. This implies that $\varphi$ may not be contracting around, say, $(x,1)$ if $\rho_3 < \theta$, $d_1$ is close to $e_3$ and $\rho_3/\theta < x < d_1$.

Starting from $z \in A^0_i$, we can apply either $f_j$ or $f_k$ to $z$, then again apply either $f_i$ or $f_k$ to $f_j(z) \in A^0_j$ and so on. Let $\sigma \in \{1,2,3\}^{\mathbb{N}}$ denote a sequence with the properties $\sigma_1 \ne i$, $\sigma_t \ne \sigma_{t-1}$, $t = 1, 2, \ldots$ (such a sequence will be called *allowed sequence*). For $z \in A^0_i$, let $f^{(t)}_\sigma(z) = f_{\sigma_t}(f_{\sigma_{t-1}}(\cdots(f_{\sigma_1}(z)\cdots))$. As the mappings $f_i$ are monotone, they are invertible and we will need, a little later, to consider mappings $f^{(-t)}_\sigma(z) = f^{(-1)}_{\sigma_t}(f^{(-1)}_{\sigma_{t-1}}(\cdots(f^{(-1)}_{\sigma_1}(z)\cdots))$.

LEMMA 4.3. *There exists a constant $\kappa > 0$ such that, for any short enough interval $[u,w] \subset {}^VA_i$, $i \in \{1,2,3\}$, and any allowed sequence $\sigma$, with the constant $\gamma \in (0,1)$ from Lemma 4.2:*

(i) $\qquad |f^{(t)}_\sigma(w) - f^{(t)}_\sigma(u)|_1 \le \gamma^t|w-u|_1 \qquad$ *for $t = 1, 2, \ldots,$*

(ii) $$e^{-\kappa|w-u|_1}\frac{|v-u|_1}{|w-u|_1} \le \frac{|f_\sigma^{(t)}(v)-f_\sigma^{(t)}(u)|_1}{|f_\sigma^{(t)}(w)-f_\sigma^{(t)}(u)|_1} \le e^{\kappa|w-u|_1}\frac{|v-u|_1}{|w-u|_1}, \qquad t=1,2,\dots,$$

*for any* $v \in (u,w)$.

PROOF. Inequalities (i) follow immediately from Lemma 4.2 since we have $f_j({}^VA_i) \subset {}^VA_j$ for all pairs $i \neq j$. For part (ii), consider any monotone function $g$ on a short interval $[a,b]$ with $|g'(x)| > \alpha_1 > 0$ and $|g''(x)| < \alpha_2$. Expanding $g$ to second order around $a$ with Taylor's theorem, we find that, for any $c \in (a,b)$,

$$\frac{g(c)-g(a)}{g(b)-g(a)} = \frac{c-a}{b-a}(1+\eta) \qquad \text{where } \eta = O(b-a)$$

and is nonzero by monotonicity of $g$. Composing this result $t$ times, we obtain

$$\frac{|f_\sigma^{(t)}(v)-f_\sigma^{(t)}(u)|_1}{|f_\sigma^{(t)}(w)-f_\sigma^{(t)}(u)|_1} = \frac{|v-u|_1}{|w-u|_1}\prod_{n=1}^{t}(1+\eta_n),$$

where $|\eta_n| < \kappa_1|w-u|_1\gamma^n$ for some $\kappa_1 > 0$ independent of $u$, $v$, $w$ and $\sigma$. This establishes the second set of inequalities. Intuitively, the idea here is that lines from points $v \in [u,w]$ to vertex $v_j$ will be almost parallel when $|w-u|_1$ is small, so relative lengths of subintervals will be about the same after mapping by $f_j$. $\square$

From Lemmas 4.2 and 4.3, it follows that there is $0 < \gamma < 1$ such that mapping $\varphi$ is uniformly contracting on $C(\gamma) := C_{\hat\imath}(\gamma) \cup C_{\hat\jmath}(\gamma) \cup C_{\hat k}(\gamma)$. Fix this $\gamma$ from now on.

LEMMA 4.4. *All points* $z$ *on a periodic orbit lie in* ${}^VA$.

PROOF. We show instead that points outside ${}^VA$, that is, $z \in J_{\hat k}$ for some $\hat k$, cannot lie on a periodic orbit because they have terminating sequences of pre-images. Suppose $J_{\hat k} \neq \varnothing$.

If neither $d_{\hat\imath}$ or $d_{\hat\jmath} \in J_{\hat k}$, then $\varphi(J_{\hat k}) \subset A^0_{\hat k} \cap {}^VA$ and no $z \in J_{\hat k}$ has a pre-image. If $d_{\hat\jmath} \in J_{\hat k}$ but $d_{\hat\imath} \notin J_{\hat k}$, then, recalling the notation $\pi(z) = x$ for $z = (x,i)$, those $z \in A^0_{\hat\imath} \cap J_{\hat k}$ with $\pi(z) < \pi(f_{\hat\imath}(d_{\hat\jmath}))$ have a single pre-image, while the other $z \in J_{\hat k}$ have none. The case with $d_{\hat\imath} \in J_{\hat k}$ but $d_{\hat\jmath} \notin J_{\hat k}$ is similar.

If both $d_{\hat\imath}$, $d_{\hat\jmath} \in J_{\hat k}$, as shown in Figure 4, then $f_{\hat\jmath}(d_{\hat\imath})$, $f_{\hat\imath}(d_{\hat\jmath}) \in J_{\hat k}$, but no $z \in \{z \in A^0_{\hat\imath} : \pi(z) > \pi(f_{\hat\imath}(d_{\hat\jmath}))\} \cup \{z \in A^0_{\hat\jmath} : \pi(z) < \pi(f_{\hat\jmath}(d_{\hat\imath}))\}$ has a pre-image under $\varphi$. Let $\sigma = \{\hat\imath, \hat\jmath, \hat\imath, \hat\jmath, \dots\}$ and observe that $f_\sigma^{(-2t)}(z) \to e_{\hat k}$ as $t \to \infty$ for any $z \in A^0_{\hat{\jmath}} \cap J_{\hat{k}}$. Hence, for $z \in A^0_{\hat{\jmath}}$ with $\pi(z) > \pi(f_j(d_{\hat{\imath}}))$, we see that $\pi(f_\sigma^{(-2t)}(z)) < \pi(f_{\hat{\jmath}}(d_{\hat{\imath}}))$ for some finite $t$ and these $z$ have only finitely many pre-images under $\varphi$. Any point in $A^0_{\hat{\jmath}} \cap J_{\hat{k}}$ can be written as $f_j(z)$ for some $z \in A^0_{\hat{\imath}} \cap J_{\hat{k}}$, so they too have finitely many pre-images in this case. We note that a forward trajectory from any $z \in J_{\hat{k}}$ either enters $^VA$ or converges onto the period two orbit on the points $(\rho_{\hat{k}}/\theta, \hat{\imath})$ and $(1 - \rho_{\hat{k}}/\theta, \hat{\jmath})$. □

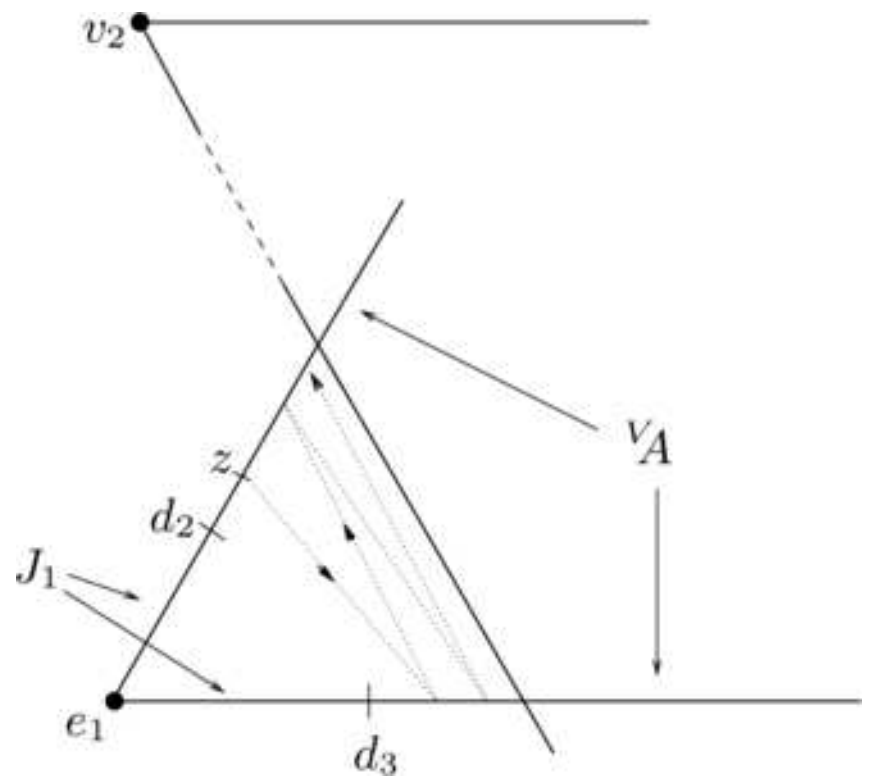

FIG. 4. *A case in the proof of Lemma 4.4.*

LEMMA 4.5. *For any starting point $z(0)$, there is $t_0 > 0$ such that $z(t) \in C(\gamma)$ for all $t \geq t_0$.*

PROOF. We see that $^VA$ is closed under $\varphi$, so once a trajectory enters it never leaves. Any trajectory that never enters $^VA$ must remain in one of the $J_{\hat{\imath}}$ and so, by Lemma 4.4, converges to the two cycle on the endpoints of that $J_{\hat{\imath}}$. In either case Lemma 4.2 implies that it enters the contracting region and remains there. □

LEMMA 4.6. *For any eventually-$m$-periodic trajectory $z(t)$, $t = 0, 1, 2, \ldots$:*

(i) *there exists an orbit $u_1, \ldots, u_m$ of period $m$ onto which the trajectory $z(t)$ converges;*
(ii) *trajectories $z'(t)$ with the same node-cycle converge onto the same orbit;*
(iii) *there can be at most one orbit having a given node-cycle.*

PROOF. (i) Note that, as $z(t)$ is *eventually-$m$*-periodic, there exists $t_0 > 0$ such that $\mathrm{Side}(z(t+m)) = \mathrm{Side}(z(t))$ for $t \geq t_0$. Now Lemma 4.5 implies there will be $t_1 > 0$ such that $z(t) \in C(\gamma)$ for all $t \geq t_1$. Since $\mathrm{Side}(z(tm+i))$ is the same for all $t \geq t_0$, Lemma 4.2 implies

$$|z((t+1)m+i) - z(tm+i)|_1 \leq \gamma^m |z(tm+i) - z((t-1)m+i)|_1$$

and, hence, $u_i = \lim_{t\to\infty} z(tm+i)$ exists. For parts (ii) and (iii), Lemmas 4.4 and 4.2 imply that, for some $r \leq m$, $|z(t) - z'(t+r)|_1 \to 0$ exponentially quickly as $t \to \infty$, so the trajectory $z'(t)$ converges onto the orbit $u_1, \dots, u_m$. The existence of another orbit $v_1, \dots, v_m$ with the same node-cycle is impossible because we can choose $z'(0) = v_1$. $\square$

Recall the notation $A^{t+1} = \varphi(A^t) = \cdots = \varphi^{(t+1)}(A^0)$. The set $A^t$ contains all points that have at least $t$ pre-images under $\varphi$. Let $\mathcal{P} = \{\varphi^{(-t)}(d_i) : i = 1, 2, 3; t = 0, 1, 2, \dots\}$ be the set of pre-images of the decision points and $\overline{\mathcal{P}}$ be its closure. Let $|z - \mathcal{P}|_1 = \inf_{z' \in \mathcal{P}} |z - z'|_1$ denote the distance from $z$ to $\mathcal{P}$—of course, $|z - \mathcal{P}|_1 = |z - \overline{\mathcal{P}}|_1$.

Remarks. If $\{d_1, d_2, d_3\} \cap A^t = \varnothing$ for some finite $t$, then $\mathcal{P}$ is finite. If $\{d_1, d_2, d_3\} \cap A^t \neq \varnothing$ for any finite $t$, then $\mathcal{P}$ can still be finite when there is a finite orbit containing one or more of the decision points. The only other possibility is that at least one decision point has infinitely many distinct pre-images and so $\mathcal{P}$ is infinite. An orbit containing a decision point must be unstable if it is of odd length $m$ and may be stable on one side if $m$ is even.

Lemma 4.7. (i) *If $\mathcal{P}$ is finite, then the triangle process has a finite number of periodic orbits and each trajectory $z(t)$, from any starting point $z(0)$, converges onto one of these orbits as $t \to \infty$.*

(ii) *Whether $\mathcal{P}$ is finite or not, if $z \notin \overline{\mathcal{P}}$ and the trajectory $z(t)$ with $z(0) = z$ converges onto the orbit $u_1, \dots, u_m$, then, for some $\varepsilon > 0$, there is a unique trajectory starting from each $z' \in \mathcal{N}_\varepsilon(z)$ and it converges onto this orbit.*

(iii) *Any orbit $u_1, \dots, u_m$ with $u_n \notin \overline{\mathcal{P}}$ for each $n = 1, \dots, m$ is stable.*

Proof. (i) As $\mathcal{P}$ is finite, we can partition $A^0$ into $|\mathcal{P}|$ open intervals $O_n$ (indexed anticlockwise from the interval containing $e_1$, say) with points of $\mathcal{P}$ as endpoints. These $O_n$ are never split by mapping with $\varphi$ as they and their images never contain decision points, so each $O_n$ is mapped into another by $\varphi$ as it is continuous and monotone except at the $d_i$. Thus, $\varphi$ induces a mapping $h$ of the indices $\{1, 2, \dots, |\mathcal{P}|\}$ into itself and, hence, $h$ has a core $K \subseteq \{1, 2, \dots, |\mathcal{P}|\}$ such that $h(K) = K$. As $h$ has no fixed points, it permutes $K$ and so factors into a product of disjoint nontrivial cycles and, hence, all trajectories starting from $z \notin \mathcal{P}$ are periodic with one of a finite number of node-cycles. Each node-cycle supports only a single orbit by Lemma 4.6 and these trajectories converge onto one of these orbits.

Now consider trajectories with $z(0) = z \in \mathcal{P}$. If, for any $t \geq 1$, we have $z(t) \in O_n$ for some $n$, then the trajectory is periodic and converges onto a finite orbit by the above argument. This only leaves trajectories with $z(t) \in \mathcal{P}$ for all $t \geq 0$. As was remarked when orbits were defined, all orbits are disjoint, so a trajectory can only remain in $\mathcal{P}$ by following a single orbit.

(ii) As $z \notin \overline{\mathcal{P}}$, the trajectory started from $z$ never hits a decision point and so never branches. As $\overline{\mathcal{P}}$ is closed, there exists $\varepsilon > 0$ such that $|z - \overline{\mathcal{P}}|_1 > \varepsilon$, so no trajectory started from any $z' \in \mathcal{N}_\varepsilon(z)$ ever branches. As $\mathcal{N}_\varepsilon(z)$ contains no pre-images of decision points, the same switching decisions are made for each of these trajectories at all times $t$, so they have the same node-cycle and by Lemma 4.6, they converge to the orbit $u_1, \dots, u_m$. Part (iii) follows immediately by considering $z(0) = u_1$. □

We are now ready to establish Theorems 3.1 and 3.2. Let $\mu$ denote Lebesgue measure on $A^0$ in the following proofs. This measure is consistent with the distance $|\cdot|_1$ we are using and any null sets in $A^0$ will remain so under change of the parameters, as the re-weighting in Lemma 4.1 is smooth.

Proof of Theorem 3.1. It suffices to show that almost all choices of decision points $d_i$ lie outside the $A^t$ they generate after some finite $t$. It follows from Lemmas 4.4 and 4.2 that $\mu(A^t) \to 0$ as $t \to \infty$ for any decision points $d_i$, but this is not sufficient since $\liminf_t A^t$ is a function of the $d_i$. What we show is that, for any given parameters $\rho_i$, there is a measure zero set of locations for the decision points where they could have infinitely many pre-images.

We start with the case $e_i \notin {}^V\!A$ or, equivalently, $J_i \neq \varnothing$, $i = 1, 2, 3$. Let $J_i^1 = f_i(J_i)$ and $J_i^{t+1} = f_i(J_j^t \cup J_k^t)$ for $i = 1, 2, 3$ and $t = 1, 2, \dots$ (here $\{j, k\} = \{1, 2, 3\} \setminus \{i\}$). These sets are well defined as $J_i^1 \subset {}^V\!A_i$ and, hence, $J_i^t \subset {}^V\!A_i$ for all $t$. Further, $J_i \cap J_i^1 = \varnothing$ and, hence, $J_i^{t+1} \cap \bigcup_{n=1}^t J_i^n = \varnothing$ for all $t$. Let $I_i^t = {}^V\!A_i \setminus \bigcup_1^t J_i^n$ and note that $I_i^{t+1} = f_i(I_j^t \cup I_k^t)$.

We initially consider only ${}^V\!A_i$. At stage $t'$, the set $\bigcup_1^{t'} J_i^n$ consists of $2^{t'+1} - 1$ disjoint intervals interleaved with the $2^{t'+1}$ intervals forming $I_i^{t'}$. The set $J_i^{t'+1}$ contains $2^{t'+1}$ intervals, each interior to one of those of $I_i^{t'}$. Recall that $\mu$ is the measure induced by the distance $|\cdot|_1$. Let $K > 0$ be such that $\mu(I_j^{t'}) < K\mu(J_j^{t'+1})$, $j = 1, 2, 3$. As $\sum_{n=1}^t \mu(J_i^n) \leq \mu({}^V\!A_i) \leq 1$, it follows that $\mu(J_i^t) \to 0$ as $t \to \infty$.

From Lemma 4.3(i), each subinterval of $I_j^{t'}$, $j = 1, 2, 3$, has length bounded by $\gamma^{t'}$ for some $\gamma < 1$. Now from Lemma 4.3(ii), it follows, by summing over all the subintervals in $I_i^{t'}$ and mapping sequences, that

$$\mu(I_i^{t'+t}) < e^{\kappa\gamma^{t'}} K\mu(J_i^{t'+t+1}), \qquad t = 2, 3, \ldots$$

and, hence, $\mu(I_i^t) \to 0$ as $t \to \infty$. The set $I_i^t$ contains all points in $A_i^0$ that could have $t$ or more pre-images, given appropriately chosen decision points. Hence, $I_i^\infty = \bigcap_1^\infty I_i^t$ contains those points that could have infinitely many pre-images and $\mu(I_i^\infty) = 0$ .

This construction is independent of the choice of the decision points and works on all three sides of ${}^V\!A$. With $I^\infty = \bigcup_1^3 I_i^\infty$, we have shown $\mu(I^\infty) = 0$, that is, $\mu$-almost all choices of the decision points have only finitely many pre-images under the mapping $\varphi$ that they define. The theorem now follows in this case from Lemma 4.7.

It remains to modify this argument in the case where one or more of the $e_i \in {}^V\!A$. For any such corner, $J_i = \varnothing$, so let $J_i^0(\varepsilon) = \{|z - e_i|_1 < \varepsilon\}$ when $e_i \in {}^V\!A$, $J_i^0 = J_i$ otherwise, where $\varepsilon > 0$ is small. The previous double mapping process will produce overlapping sets, so we modify it by setting $J_i^1 = f_i(J_i^0) \setminus (J_j^0 \cup J_k^0)$ and then $J_i^{t+1} = f_i(J_j^t \cup J_k^t) \setminus (J_j^0 \cup J_k^0)$ for $i = 1,\ 2,\ 3$ and $t = 1,\ 2,\ \ldots$. This change ensures the previous disjointness properties and the construction goes through without further change, except that now some of the sub-intervals of $I_i^t$ may be empty. As before, we conclude that $\mu(I^\infty) = 0$, though here $I^\infty$ does depend upon $\varepsilon$.

To relate this construction to the existence of pre-images, let $\mathcal{D}_n = \{(d_1, d_2, d_3) :$ all $z \in J_i^0(1/n), i = 1, 2, 3$, have no legitimate pre-images under $\varphi\}$. Then $\bigcup_n \mathcal{D}_n = A_1^0 \times A_2^0 \times A_3^0$. For any fixed $n$, $\mathcal{D}_n^\infty = \mathcal{D}_n \cap I^\infty$ has $\mu(\mathcal{D}_n^\infty) = 0$ and hence $\mu(\bigcup_1^\infty \mathcal{D}_n^\infty) \le \sum_1^\infty \mu(\mathcal{D}_n^\infty) = 0$. Now the result follows as before. $\square$

Proof of Theorem 3.2. We consider first the case where $\mathcal{P} = \{\varphi^{(-t)}(d_i) : i = 1, 2, 3; t = 0, 1, 2, \ldots\}$, the set of pre-images of the decision points, is finite. A point $z \in \mathcal{P}$ is of *type* $i$ when it is a pre-image of $d_i$. $\mathcal{P}$ splits $A^0$ into closed intervals $M_n$, $n = 1,\ 2,\ \ldots, |\mathcal{P}|$ which are never split by iterated mappings with $\varphi$. Hence, each $M_n$ contains points from at most one periodic orbit and as $\varphi$ is contracting on ${}^V\!A$, each $M_n$ can contain only a single point from an orbit—in cases where the shared endpoint of $M_n$ and $M_{n+1}$ lies on an orbit, that is, this orbit contains one or more decision points, there can be at most one other orbit with a point in $M_n \cup M_{n+1}$.

The number of orbits is limited by each $d_i$ being at the boundary of just two intervals. Consider an orbit with a point $u \in M_t$ which has an endpoint $\varphi^{(-r)}(d_i)$. The point $\varphi^{(r)}(u)$ lies on the same orbit but occupies one of the $M_n$ neighboring $d_i$, or equals $d_i$. By the preceding paragraph, there are at most six orbits.

We can reduce the bound as there are at least two intervals $M_n$ with endpoints of different type, $ij$ and $ik$, say. Suppose the first is $M_t$ with endpoints $\varphi^{(-r)}(d_i)$ and $\varphi^{(-s)}(d_j)$ with $r < s$. The interval $\varphi^{(r)}(M_t)$ has endpoints $d_i$ and $\varphi^{(r-s)}(d_j)$ so it neighbors $d_i$. $\varphi^{(s)}(M_t)$ sits inside an interval neighboring $d_j$. If $M_t$ contains a point from an orbit, then this orbit occupies two of the six intervals neighboring the decision points. If $M_t$ contains no point from any orbit, then neither can one of the intervals neighboring $d_i$. Either way, the maximum number of orbits is now only five. Repeat the argument with the interval with endpoints $ik$ and there are at most four orbits in the stable case, that is, where $\mathcal{P}$ is finite.

We now extend this argument to the case where $\mathcal{P}$ is infinite—let $\mathcal{P}' = \{z \in \mathcal{P} : z \text{ has infinitely many pre-images}\}$. Suppose there are some finite orbits and consider any of them. It cannot include any points in $\mathcal{P}'$, so choose $\varepsilon > 0$ small enough that $|u_n - d_i|_1 > \varepsilon$ for every $u_n$ on the orbit and each $d_i \in \mathcal{P}'$. It follows that $|u_n - \varphi^{(-r)}(d_i)|_1 > \varepsilon$ for every $u_n$ on the orbit and every point in $\mathcal{P}'$ for suppose this was not true. Select the $\varphi^{(-r)}(d_i)$ with minimal $r$ within $\varepsilon$ of $u_n$. By Lemma 4.3(i), it follows that

$$|\varphi^{(r)}(u_n) - d_i|_1 \le \gamma^r |u_n - \varphi^{(-r)}(d_i)|_1 < \varepsilon,$$

in contradiction to our choice of $\varepsilon$. Now our argument for the stable case applies with the added possibility that one or more $d_i$ may be limit points (from one or both sides) of $\mathcal{P}$ further restricting the opportunities for orbits. □

We show later that it is possible for finite orbits to exist simultaneously with nonperiodic trajectories in this nonstable case. In our study of this model we simulated the triangle process over the whole range of its parameters $\rho_i$ and $d_i$. To avoid problems with rounding, we used an algebraic representation of the process which is described below in Section 5. The representation there of the decision points has to be finite in practice, so the results of the computer analysis were sometimes inconclusive, for example, the implemented algorithm sometimes produced trajectories that were not trapped by an orbit in the number of steps we could accurately calculate. Otherwise, the algorithm identified either 1, 2 or 3 orbits and never more. Therefore, we *conjecture* that Theorem 3.1 can be strengthened, by replacing *at most four orbits* by *three orbits*. We did not manage to prove this analytically, though. That three orbits can exist is easily seen after a little numerical work with Figure 5 as a guide.

**5. The triangle process in the nonstable case.** We now show that there are locations for the decision points where at least one of them has infinitely many pre-images and nonperiodic trajectories exist. This argument necessarily uses yet another way of describing the triangle process. We will consider only the case where each $J_i = \varnothing$, but our construction can be carried out in other cases too, with some limitations.

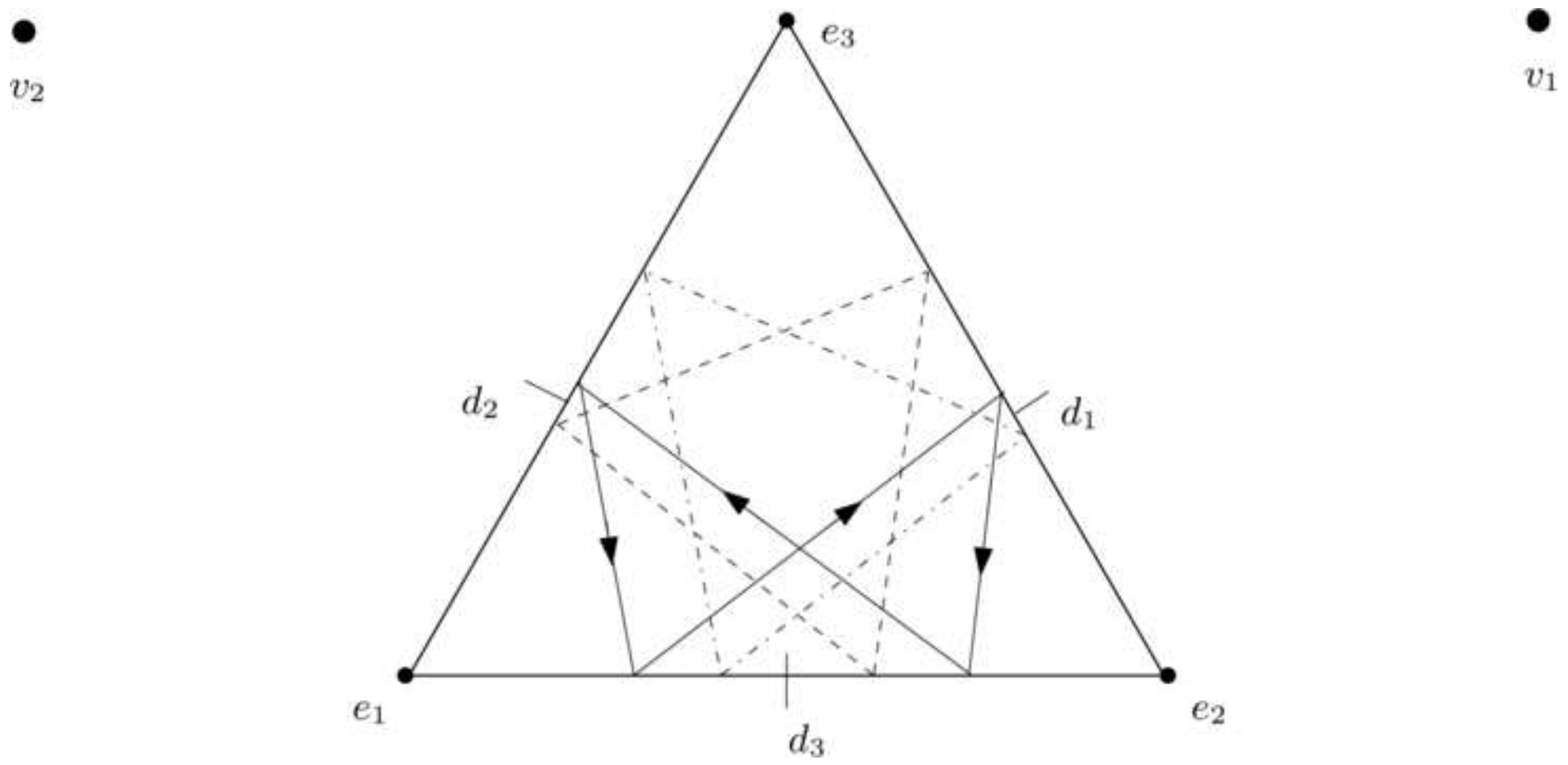

FIG. 5. *An example with three orbits.*

We first construct a binary representation $z = \hat{\imath} : x_1 x_2 x_3 \ldots$ with each $x_t \in \{0, 1\}$ for all points of $A^0$. For each $z \in A^0_{\hat{\imath}}, \hat{\imath} = 1, 2, 3$, set [recall the notation $\pi(\cdot)$ from Section 3]

$$x_1 = \begin{cases} 0, & \text{if } \pi(z) < \pi(f_{\hat{\imath}}(e_{\hat{\imath}})), \\ 1, & \text{if } \pi(z) \geq \pi(f_{\hat{\imath}}(e_{\hat{\imath}})). \end{cases}$$

To continue this construction, set $B^1_{\hat{\imath}} = \{f_{\hat{\imath}}(e_{\hat{\imath}})\}$ and $B^{t+1}_{\hat{\imath}} = f_{\hat{\imath}}(B^t_{\hat{\jmath}} \cup B^t_{\hat{k}})$ for $\hat{\imath} = 1,\ 2,\ 3$ and $t = 1,\ 2,\ \ldots.$

Each set $B^t_{\hat{\imath}}$ contains $2^{t-1}$ points which, by monotonicity of the mappings $f_{\hat{\imath}}$, interleave the $1 + 2^{t-1}$ points of $\bigcup_{r=1}^{t-1} B^r_{\hat{\imath}}$ and, hence, split in two each of the $2^{t-1}$ intervals created up to stage $t-1$. Label the $b_r \in B^t_{\hat{\jmath}} \cup B^t_{\hat{k}}$ clockwise and the $a_r \in \bigcup_{r=1}^{t-1} B^r_{\hat{\imath}}$ anticlockwise so that $a_0 = e_j$ and

$$\pi(a_{r-1}) < \pi(f_{\hat{\imath}}(b_r)) < \pi(a_r), \qquad r = 1, 2, \ldots, 2^{t-1},$$

and for $z \in [a_{r-1}, a_r)$, set

$$x_t = \begin{cases} 0, & \text{if } \pi(z) < \pi(f_{\hat{\imath}}(b_r)), \\ 1, & \text{if } \pi(z) \geq \pi(f_{\hat{\imath}}(b_r)). \end{cases}$$

This encodes every point of $A^0$ with a unique binary code up to the usual indeterminacy for points terminating with infinite strings of zeros or ones, for example, $e_1 = 2 : 111 \cdots = 3 : 000 \cdots$.

This binary encoding can be used to define a distance between points $z = i : x_1 x_2 \cdots, z' = i : x'_1 x'_2 \cdots$ by

$$|z - z'|_b \equiv \sum_{t=1}^{\infty} |x_t - x'_t| \cdot 2^{-t} \tag{11}$$

and we note the following: (i) ordering of points on $A_i^0$ by $\pi(z)$ corresponds to lexicographic ordering on the binary sequence $x$; (ii) any sequence of points which is convergent under $|\cdot|_1$ remains convergent under $|\cdot|_b$ and vice versa; (iii) this metric induces a measure $\mu_b$ on $A_i^0$.

In defining the required decision points, we will make much use of the inverse mapping $\psi$, where

$$\psi(z) = f_i^{(-1)}(z), \qquad z \in A_i^0$$

for $i = 1, 2, 3$. When $\varphi^{(-1)}(z)$ exists, then we will say $\psi$ is *legitimate* and we have $\psi(z) = \varphi^{(-1)}(z)$. The mapping $\psi$ can be described quite simply using the binary encoding. For $z = \hat{\imath} : x_1 x_2 \cdots$ with $x_1 = 0$, then $\psi(z) = \hat{k} : y_1 y_2 \cdots$, where $y_1 = 1 - x_2$, $y_2 = 1 - x_3$ and so on. In general, let $\bar{w} = 1 - w$ for $w \in \{0, 1\}$, $\bar{x} = \bar{x}_1 \bar{x}_2 \cdots$ and with this notation $\psi$ maps

$$\hat{\imath} : 0x \longrightarrow \hat{k} : \bar{x} \quad \text{and} \quad \hat{\imath} : 1x \longrightarrow \hat{\jmath} : \bar{x} \tag{12}$$

for $\hat{\imath} = 1,\ 2,\ 3$.

The forward mapping $\varphi$ thus can be represented as

$$\hat{\imath} : x \longrightarrow \hat{\jmath} : 0\bar{x} \qquad \text{if } x \leq d_{\hat{\imath}} \quad \text{and} \quad \hat{\imath} : x \longrightarrow \hat{k} : 1\bar{x} \qquad \text{if } x \geq d_{\hat{\imath}}.$$

Now represent each point of $A^0$ by a point in the interval $[0, 1]$ using the mapping

$$i : x_1 x_2 \cdots \longrightarrow \frac{i-1}{3} + \frac{1}{3} \sum_{j=1}^{\infty} 2^{-j} x_j$$

and ignore, for the moment, the branching at the decision points. The forward mapping $\varphi$ becomes

$$\varphi(x) = \begin{cases} -\frac{1}{2}x + 1, & \text{for } x \in [\tilde{d}_1, \tilde{d}_2) \cup [\tilde{d}_3, 1), \\ -\frac{1}{2}x + \frac{1}{2}, & \text{for } x \in [0, \tilde{d}_1) \cup [\tilde{d}_2, \tilde{d}_3), \end{cases} \tag{13}$$

for $x \in [0, 1)$, where $\tilde{d}_1$, $\tilde{d}_2$ , $\tilde{d}_3 \in [0, 1)$ are the points corresponding to the decision points. The mapping above can be regarded as a particular case of the affine interval exchange transformation model (see, e.g., [2, 7]). The mapping defined by (13) is contracting and is not order-preserving, which makes it different from the usual and well-studied case of interval exchange transformations (see, e.g., [3, 11] and references therein). Note also that, for interval exchange transformations with contracting and/or flips, it is natural that cycles do exist (see [3]), as happens in our model.

5.1. *A decision point with infinitely many pre-images.* We are now ready to construct decision points providing a nonstable case. We will use $q = 1001$ and $r = 0110$, two four digit sequences such that $\bar{q} = r$ and such that $\psi$ successively maps $z = \hat{\imath}: qx$ and $z = \hat{\imath}: rx$ as follows:

$$z = \hat{\imath}: 1001x \longrightarrow \hat{\jmath}: 110\bar{x} \longrightarrow \hat{k}: 01x \longrightarrow \hat{\jmath}: 0\bar{x} \longrightarrow \hat{\imath}: x = \psi^{(4)}(z)$$

and, similarly,

$$z = \hat{\imath}: 0110x \longrightarrow \hat{k}: 001\bar{x} \longrightarrow \hat{\jmath}: 10x \longrightarrow \hat{k}: 1\bar{x} \longrightarrow \hat{\imath}: x = \psi^{(4)}(z).$$

The decision points have the form

$$\begin{aligned} d_1 &= 1: qrx\cdots, \\ d_2 &= 2: 1010100000\cdots, \\ d_3 &= 3: 0100000\cdots, \end{aligned} \tag{14}$$

where $x$, the infinite tail of the sequence for $d_1$, will be constructed recursively at a later point to ensure that $d_1$ has infinitely many legitimate $\psi$-images. As the notation $d_i$ identifies the triangle side, we will sometimes use it to denote just the binary sequence. Note that, under lexicographic ordering, $r = 0110 < 1001 = q$.

For $d_1$ to have infinitely many legitimate $\psi$-images, it must certainly have one and the intervals with at least one legitimate $\psi$-image are

$$\begin{aligned} &\text{on } A_1^0 \qquad [0\overline{d}_3, 1\overline{d}_2] = [0101111\cdots, 1010101111\cdots], \\ &\text{on } A_2^0 \qquad [0\overline{d}_1, 1\overline{d}_3] = [0rq\bar{x}\cdots, 110111111\cdots], \\ &\text{on } A_3^0 \qquad [0\overline{d}_2, 1\overline{d}_1] = [00101011111\cdots, 1rq\bar{x}\cdots], \end{aligned}$$

which the reader can readily check certainly contain $d_1$ and the $\psi^{(t)}(d_1)$ for $t = 1, 2, \ldots, 8$. We will construct the rest of the binary sequence for $d_1$ so that its $\psi$-images always fall in these three intervals by finding a sequence of quadruples $q$ and $r$ that guarantees legitimacy.

*Extended legitimacy property.* The sequence $d_1 = y_1y_2y_3\cdots$ where $y_1 = q$, $y_2 = r$ and $y_t \in \{q, r\}$ for $t = 3, 4, \ldots$ has *extended legitimacy* if, under lexicographical ordering,

(a) when $y_t = r$, then $y_{t+1}y_{t+2}\cdots > d_1$,
(b) when $y_t = q$, then $y_{t+1}y_{t+2}\cdots < d_1$,

where $q > r$ as remarked above.

At each $t \geq 3$ this property ensures the legitimacy of $\psi^{(4t-r)}(d_1)$ for $r = 3$, 2, 1, 0. We will now construct an aperiodic sequence with extended legitimacy and show that this means there are decision points leading to the existence of aperiodic orbits which is crucial for establishing Theorem 3.3.

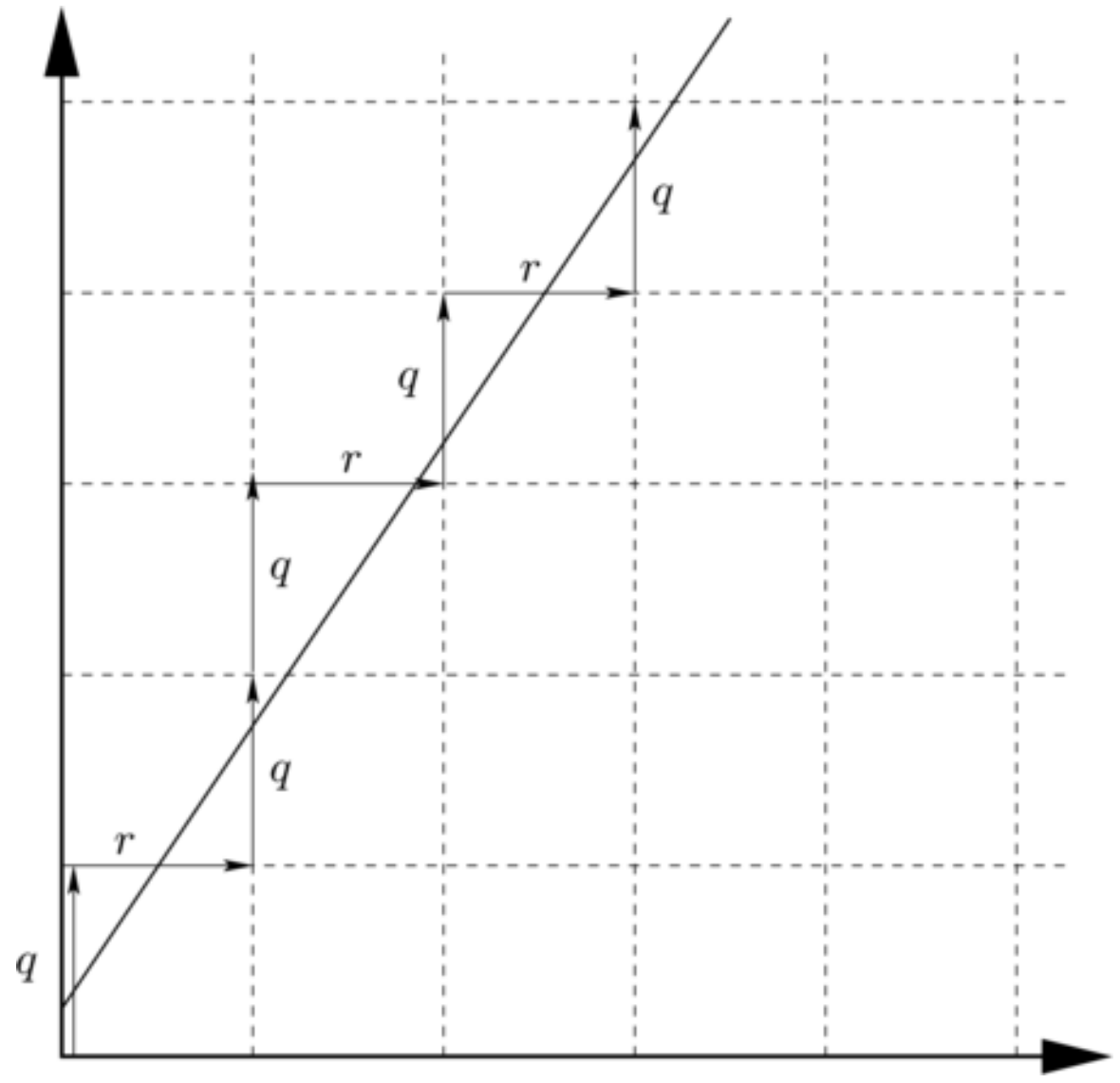

Fig. 6. *Staircase diagram for sequence qrqqrqrq....*

Fix an *irrational* $\alpha \in (1,2)$ and some $\beta \in [-1,1]$. Let $y = y(\alpha,\beta) = y_1y_2y_3\cdots$ be a sequence of quadruples $q$ and $r$ with $y_1 = q$, $y_2 = r$ and for $t = 3,4,\ldots,$

$$y_{t+1} = \begin{cases} q, & \text{if } 1 + Q_t < \alpha(1+R_t) + \beta, \\ r, & \text{otherwise,} \end{cases}$$

where $Q_t = Q_t(y) = \sum_{r=1}^{t} \mathbf{1}_{\{y_r=q\}}$ and $R_t = t - Q_t$, that is, $Q_t$ is the number of times $q$ appears in the sequence $y_1y_2\cdots y_t$. This sequence is defined using successive rational approximations to the irrational number $\alpha$ as shown in Figure 6, which shows a *staircase* approximating the line $y = \alpha x + \beta$.

Lemma 5.1. *There are uncountably many sets of decision points where at least one of them has infinitely many pre-images.*

Proof. We first show that, for any irrational $\alpha > 0$, the sequence $y' = y(\alpha,\beta')$ is lexicographically less than $y'' = y(\alpha,\beta'')$ whenever $-1 \le \beta' < \beta'' \le 1$.

Suppose $y'$ and $y''$ are not equal. Then there is $t \ge 1$ such that $y'_n = y''_n$ for $n \le t$, but $y'_{t+1} \ne y''_{t+1}$. Evidently, $Q_t(y') = Q_t(y'')$ and $R_t(y') = R_t(y'')$, so the only way for $y'_{t+1} \ne y''_{t+1}$ is

$$\begin{aligned} \alpha(1+R_t(y')) + \beta' &< 1 + Q_t(y') \\ &= 1 + Q_t(y'') \\ &< \alpha(1+R_t(y'')) + \beta'', \end{aligned}$$

in which case $y'_{t+1} = r < q = y''_{t+1}$. As the fractional parts of $m\alpha$, $m = 1, 2, \dots,$ are dense in $(0, 1)$, the inequalities $\beta' < n - m\alpha < \beta''$ are satisfied for infinitely many pairs of integers $m$, $n$ and so the above inequalities also show that $y' = y''$ is not possible.

Now consider the sequence $y = y(\alpha, 0)$ which starts $qr \cdots$. Each $r$ is followed by a $q$ and there cannot be more than two $q$s in succession as $\alpha < 2$. Suppose $y_t = r$. Then the sequence $y_{t+1}y_{t+2}y_{t+3}\cdots$ coincides with $y(\alpha, \beta_t)$, where

$$\beta_t = \alpha R_t - Q_t > 0$$

and so $y_{t+1}y_{t+2}y_{t+3}\cdots > y$. If instead $y_t = q$, then either $y_{t+1} = r$ so that immediately we have $y_{t+1}y_{t+2}y_{t+3}\cdots < y$ or $y_{t+1} = q$, in which case $y_{t+1}y_{t+2}y_{t+3}\cdots$ coincides with $y(\alpha, \beta_t)$, where

$$\beta_t = \alpha R_t - Q_t < 0$$

and so $y_{t+1}y_{t+2}y_{t+3}\cdots < y$. In all cases we see that the sequence $y$ has extended legitimacy at $t$ for all $t \geq 3$. As $\alpha$ is irrational, $\alpha R_t - Q_t = 0$ is not possible and, hence, there is no $t$ such that $y_{t+1}y_{t+2}y_{t+3}\cdots = y$ from which it follows immediately that $y$ is aperiodic.

To establish this lemma, let $d_1 = 1 : y(\alpha, 0)$ and $d_2$, $d_3$ be as defined in (14) for some irrational $\alpha \in (1, 2)$. By construction, $d_1$ has an infinite sequence of legitimate $\psi$-images or, equivalently, an infinite sequence of pre-images under $\varphi$. □

REMARKS. For the above set of the decision points, there is also a period 3 orbit: $(a, c, b)$, with the points of the orbit given by

$$\begin{aligned} &1 : 101010 \cdots 10 \cdots, \\ &3 : 101010 \cdots 10 \cdots, \\ &2 : 101010 \cdots 10 \cdots. \end{aligned} \tag{15}$$

We mention here that originally we were not sure that the nonstable case was possible and unsuccessfully tried to prove this using the triangle process description of this section. The specific combination of decision points producing the nonperiodic orbits described above was discovered by studying computer output (of a programme written to generate sample trajectories) to assess various lines of investigation.

5.2. *The nonstable case has $\mu_b$-measure* 0. The measure $\mu_b$ on $A^0$ induced by the metric $|\cdot|_b$ is generally singular with respect to Lebesgue measure, but the set of decision points with infinitely many $\varphi$ pre-images has $\mu_b$-measure 0.

Suppose we select the decision points uniformly, that is, each digit in the binary sequence is 0 or 1 with chance 1/2 independently of the other digits. Then $\mu_b$ a.s. we can write, for some $k \geq 0$,

$$d_2 = 2 : \underbrace{11 \cdots 11}_{k \text{ times}} 0 \cdots .$$

This means the interval in $A_1^0$ with legitimate $\psi$-images has left endpoint

$$\underbrace{00 \cdots 00}_{k+1 \text{ times}} 1 \cdots .$$

The Borel–Cantelli lemma tells us that $\mu_b$ a.s. there is a subsequence in the binary sequence for $d_1$ of the form

$$y = \underbrace{00 \cdots 00}_{m \text{ times}} \; 10 \; \underbrace{00 \cdots 00}_{m \text{ times}} \; 10 \; \underbrace{00 \cdots 00}_{m \text{ times}},$$

where $m > k+3$ is even. Now consider the sequence $\psi^{(t)}(d_1)$, $t \geq 1$. For some $r$, $\psi^{(r)}(d_1)$ starts with either the above sequence $y$ or with $\bar{y}$ and for some $s \leq 2m+6$, $\psi^{(t+s)}(d_1)$ is a point on $A_1^0$ with binary sequence starting with at least $k+2$ zeros and so its next $\psi$-image is not legitimate. It follows immediately that the set of decision points $d_i$ with infinitely many pre-images has $\mu_b$-measure 0.

5.3. *Properties of the deterministic system when $\mathcal{P}$ is infinite.* For the triple of the decision points constructed in (14), $d_1$ has infinitely many pre-images and since $\psi^{-5}(d_2) = 1:000\cdots = e_2$ and $\psi^{-2}(d_3) = 1:111\cdots = e_3$, both $d_2$ and $d_3$ have finitely many pre-images under $\varphi$, unlike $d_1$. Throughout this section, the triple of the decision points is assumed to be such that $d_1$ has infinitely many pre-images, while $d_2$ and $d_3$ have only finitely many. We also assume in this subsection that ${}^{V}\!A = A^0$ so that $\varphi$ is contracting everywhere.

Let $\mathcal{P}_1 = \{d_1, \psi(d_1), \psi^{(2)}(d_1), \ldots\}$ be the infinite set of pre-images of $d_1$ under $\varphi$ which, by assumption, are all legitimate. We will study the properties of this set. Let $\overline{\mathcal{P}}_1$ denote the closure of the set $\mathcal{P}_1$.

LEMMA 5.2. *Each $u \in \mathcal{P}_1$ is a limiting (accumulation) point of $\mathcal{P}_1$.*

PROOF. It is sufficient to show that $d_1$ is a limiting point for $\mathcal{P}_1$. Indeed, since the mapping $\psi$ doubles the distance $|\cdot|_b$ between points, if $d_1$ is limiting for a sequence $u_1, u_2, \ldots \subset \mathcal{P}_1$, then $\psi^{(i)}(d_1)$ is a limiting point for $\psi^{(i)}(u_1)$, $\psi^{(i)}(u_2), \ldots \subset \mathcal{P}_1$ [as mentioned just after (11), the distances $|\cdot|_1$ and $|\cdot|_b$ are topologically equivalent so sets of points generate the same limit points under both].

Consider any point $u \in A^0$ which is a limiting point for $\mathcal{P}_1$, that is, there is a sequence of indices $t_1 < t_2 < t_3 < \cdots$ such that $\lim_n \psi^{(t_n)}(d_1) = u$. No such $u$ can be a pre-image of $d_2$ or $d_3$, as for each $t$, $A^0 \setminus A^t$, the set of points with fewer than $t$ legitimate pre-images, is open. Hence, around each pre-image of $d_2$ or $d_3$ there is an open interval of points which have only finitely many pre-images and so cannot contain any $\psi^{(t)}(d_1)$.

If $u \in \mathcal{P}_1$, then $u = \psi^{(m)}(d_1)$ for some $m$ which is unique (otherwise $u$ is in a finite orbit and $\mathcal{P}_1$ cannot be infinite). Without loss of generality, we can assume $m < t_1$ and then, since the mapping $\varphi^{(m)}$ is continuous except at $\psi^{(n)}(d_1)$ for $n < m$, $d_1 = \phi^{(m)}(u) = \lim_n \psi^{(t_n - m)}(d_1)$.

If $u \notin \mathcal{P}_1$, fix a small $\varepsilon > 0$ and consider the $\varepsilon$-neighborhood $\mathcal{N}_\varepsilon(u)$ of $u$ as defined in (8). Assume that $\varepsilon$ is so small that there are no pre-images of $d_2$ or $d_3$ in $\mathcal{N}_\varepsilon(u)$. Let $m$ be the smallest index such that $\psi^{(m)}(d_1) \in \mathcal{N}_\varepsilon(u)$. Then $\mathcal{N}_\varepsilon(u)$ contains infinitely many points of $\mathcal{P}_1$ and each is a pre-image of $d_1$ with the index of at least $m$. Now map $\mathcal{N}_\varepsilon(u)$ with $\varphi^{(m)}(\cdot)$ which sends $\psi^{(m)}(d_1)$ to $d_1$ and all the points of the neighborhood $\mathcal{N}_\varepsilon(u)$ follow the same trajectory. As $\varphi^{(m)}(\cdot)$ is contracting here, $\varphi^{(m)}(\mathcal{N}_\varepsilon(u))$ is an interval within $\mathcal{N}_\varepsilon(d_1)$. This implies that there are infinitely many points of $\mathcal{P}_1$ in $\mathcal{N}_\varepsilon(d_1)$. Letting $\varepsilon \downarrow 0$, we obtain that $d_1$ is a limiting point for $\mathcal{P}_1$. $\square$

The following statement immediately follows from Lemma 5.2, properties of perfect sets and the fact that $\mathcal{P}_1$ is infinite and countable.

COROLLARY 5.1. *The set $\overline{\mathcal{P}}_1$ is a perfect set ($\equiv$ it is closed and every point of it is an accumulation point). Therefore, $\overline{\mathcal{P}}_1$ is uncountable, and hence, $\overline{\mathcal{P}}_1 \setminus \mathcal{P}_1$ is also uncountable.*

Recall that $\mathcal{P}$ is the union of $\mathcal{P}_1$ and the set of $d_2$ and $d_3$ and their finitely many pre-images. Then $\overline{\mathcal{P}}$ is the union of $\overline{\mathcal{P}}_1$ and the set of $d_2$ and $d_3$ and their pre-images, and is also closed.

LEMMA 5.3. *If $u \in \overline{\mathcal{P}}_1$, then $u$ does not belong to any finite orbit and no trajectory $z(t)$ with $z(0) = u$ can be periodic.*

PROOF. If $u \in \mathcal{P}_1$, then $\varphi^{(n)}(u) = d_1$ for some $n$. However, if $u$ belonged to a finite orbit, then $\varphi^{(m)}(u) = u$ for some $m$ and, hence, $\varphi^{(m)}(d_1) = d_1$, contradicting the fact that $d_1$ has more than $m$ *distinct* pre-images under $\psi_r = \varphi^{(-1)}$.

Now suppose $u \in \overline{\mathcal{P}}_1 \setminus \mathcal{P}_1$ and $u$ belongs to a finite orbit. As $d_2, d_3 \notin A^n$ for some finite $n$ and $A^n$ is closed, we can choose $\varepsilon > 0$ so that $|d_i - A^n|_1 > \varepsilon$ for $i = 2, 3$. Since $\overline{\mathcal{P}}_1 \subset A^n$, it follows that the neighborhood $\mathcal{N}_\varepsilon(u)$ contains no pre-images of $d_2$ or $d_3$. Let $n = n(\varepsilon)$ be the smallest index such that $\psi^{(n)}(d_1) \in \mathcal{N}_\varepsilon(u)$. As $\varphi^{(n)}$ does not split this neighborhood, $\varphi^{(n)}(\mathcal{N}_\varepsilon(u)) \subset$ $\mathcal{N}_{2\varepsilon}(d_1)$ and so there is a point of the orbit, $u'$, say, with $|u'-d_1|_1 < 2\varepsilon$. Letting $\varepsilon \downarrow 0$, we obtain a contradiction, since this orbit contains only $m$ distinct points and does not contain $d_1$.

Finally, consider any trajectory $z(t)$ with $z(0)=u$ and suppose it is periodic. By Lemma 4.6(i), $z(t)$ converges onto an orbit $w_1, \ldots, w_m$ and we have just shown there exists $\varepsilon > 0$ such that $\mathcal{N}_\varepsilon(w_1) \cap \overline{\mathcal{P}}_1 = \varnothing$. However, $z(t) = \varphi^{(t)}(u) \in \mathcal{N}_\varepsilon(w_1)$ for some large $t$, while simultaneously, by Lemma 5.2, $\varphi^{(t)}(u)$ is a limiting point and, hence, an element of $\overline{\mathcal{P}}_1$ which is impossible. □

Proof of Theorem 3.3. Lemma 5.1 explains how to construct $d_1$ with an infinite sequence of pre-images under $\varphi$. Also, by Lemma 5.1, the binary sequence for $d_1$ is not periodic and, hence, by Lemma 5.3, there are nonperiodic trajectories $\{z(t): t \geq 0\}$ for the triangle process. □

This completes the main task of this section, but it is possible to say more about the mixture of periodic and nonperiodic trajectories. We know from Lemma 4.6(ii) that all periodic trajectories avoiding decision points are stable. It is easy to see that *there can be* stable periodic trajectories even if $\mathcal{P}$ is infinite. For example, for the decision points given by (14), the sequence in (15) is an orbit with period 3.

We say that a point $u \in A_i^0$ is an $R$-limit ($L$-limit resp.) for $\mathcal{P}_1$ if there is a sequence of points $u_1, u_2, \ldots$ in $A_i^0 \cap \mathcal{P}_1$ such that $\pi(u_n) \downarrow u$ [$\pi(u_n) \uparrow u$ resp.], that is, for an $R$-limit, $u_{n+1}$ is clockwise of $u_n$ for each $n$. The next result strengthens Lemma 5.2.

Lemma 5.4. *Each point $u \in \mathcal{P}_1$ is limiting for $\mathcal{P}_1$ on both sides, that is, it is both an $R$-limit and an $L$-limit for $\mathcal{P}_1$.*

Proof. It suffices to show that $d_1$ is both an $R$-limit and an $L$-limit, since all other points of $\mathcal{P}_1$ are $\psi$-images of $d_1$.

By Lemma 5.2, $d_1$ is a limiting point for $\mathcal{P}_1$, but suppose that it is just an $R$-limit and not an $L$-limit. Then there is $\delta > 0$ such that:

(i) there are no $\mathcal{P}_1$ points on the interval $(d_1-\delta, d_1)$, and
(ii) there are no pre-images of $d_2$ or $d_3$ in the segment $[d_1, d_1+\delta]$.

Since mapping $\psi$ reverses the orientation of intervals, if $u = \psi^{(n)}(d_1) \in A_1^0$ and $n$ is even, then $u$ is an $R$-limit and not an $L$-limit, while if $n$ is odd, then $u$ is an $L$-limit and not an $R$-limit.

Fix $0 < \delta_1 < \delta$. Denote $d_1^n := \psi^{(n)}(d_1), n = 1, 2, \ldots$, and define

$$r := \min\{n : \psi^{(n)}(d_1) \in (d_1, d_1+\delta_1]\}, \tag{16}$$

so $d_1^r$ is the point of $\mathcal{P}_1 \cap [d_1, d_1 + \delta_1]$ with the smallest index. If $r$ is even, then we obtain a contradiction since $\varphi^{(r)}((d_1, d_1^r])$ is not split during these mappings and $\varphi^{(r)}(d_1^r) = d_1$. Since $r$ is even, $\varphi^{(r)}$ preserves orientation and is contracting so that $\varphi^{(r)}((d_1, d_1^r]) \subset (d_1 - \delta, d_1]$, which contradicts the assumption that $(d_1 - \delta, d_1) \cap \mathcal{P}_1 = \varnothing$.

Therefore, in any small interval $(d_1, d_1 + \delta_1)$, $d_1^r$ has $r$ odd and there will be infinitely many $d_1^n$ with odd $n$ within $(d_1, d_1 + \delta_1)$. This implies there will also be $d_1^m$ with *even* $m$ in $(d_1, d_1 + \delta_1)$ since $\varphi^{(r)}$ reverses orientation and is contracting, so $\varphi^{(r)}((d_1, d_1^r]) \subset [d_1, d_1^r)$ and all the $d_1^n \in (d_1, d_1^r]$ with odd $n > r$ have been mapped to $d_1^{n-r}$ with $n - r$ even.

Now choose $d_1^m$ in $(d_1, d_1 + \delta)$ with the smallest possible even index $m \geq 2$ and choose odd index $t < m$ such that $d_1^t \in (d_1, d_1^m)$ with $d_1^t$ the closest of such points to $d_1^m$. Such a $t \geq r$ exists by the arguments following (16) above. The interval $(d_1^t, d_1^m]$ contains no $d_1^n$ with index less than $m$ so, as $\varphi^{(t)}$ reverses orientation and is contracting, $\varphi^{(t)}(d_1^m) = d_1^{m-t} \in (d_1 - \delta, d_1)$, contradicting the assumption that $(d_1 - \delta, d_1) \cap \mathcal{P}_1 = \varnothing$. $\square$

LEMMA 5.5. *Fix a point $u$ and suppose $u \notin \overline{\mathcal{P}}$. Then $u$ belongs to the segment $[\tilde{l}(u), \tilde{r}(u)] = [\tilde{l}, \tilde{r}]$ such that each of $\tilde{l}$ and $\tilde{r}$ is either a pre-image of $d_2$ or $d_3$, or belongs to $\overline{\mathcal{P}}_1$ and $\tilde{l}(u)$ is the largest of such points and $\tilde{r}(u)$ is the smallest of such points:*

$$\tilde{l}(u) = \max\{\pi(s) : \mathrm{Side}(s) = \mathrm{Side}(u),\ \pi(s) < \pi(u),\ s \in \overline{\mathcal{P}}\},$$

$$\tilde{r}(u) = \min\{\pi(s) : \mathrm{Side}(s) = \mathrm{Side}(u),\ \pi(s) > \pi(u),\ s \in \overline{\mathcal{P}}\}.$$

*Moreover, the trajectory started from $u$ can be periodic only if both $\tilde{l}(u)$ and $\tilde{r}(u)$ are pre-images of $d_2$ or $d_3$.*

PROOF. The existence of $\tilde{l}(u)$ and $\tilde{r}(u)$ follows from the fact that $\overline{\mathcal{P}}$ is closed. Also, $\pi(\tilde{l}(u)) < \pi(\tilde{r}(u))$ since the complement of the set $\overline{\mathcal{P}}$ is open.

Suppose, for example, that $\tilde{l}(u)$ is not a pre-image of $d_2$ or $d_3$, that is, $\tilde{l}(u) \in \overline{\mathcal{P}}_1$. Since there are no pre-images of the decision points on $[\tilde{l}(u), u]$, $u$ and $\tilde{l}(u)$ can follow the same trajectory under mapping $\varphi$. On the other hand, by Lemma 5.3, no trajectory started from $\tilde{l}(u)$ is periodic, hence, the trajectory starting from $u$ is not periodic. $\square$

Lemma 5.5 implies that the set of points from which periodic trajectories can start—call this the *periodic set*—consists of finitely many intervals, since $d_2$ and $d_3$ have only finitely many pre-images. Also, the periodic set has measure smaller than the total length of the sides of the triangle, since $\mathcal{P}_1$ is infinite. At the same time, Lemma 5.4 yields that, for any $u \notin \overline{\mathcal{P}}$, if $\tilde{l}(u) \in \overline{\mathcal{P}}_1$, then $\tilde{l}(u) \notin \mathcal{P}_1$, and the same is true for $\tilde{r}(u)$.

Call the set of points $u$ for which one of $\tilde{l}(u)$ or $\tilde{r}(u)$ is in $\overline{\mathcal{P}}_1$, while the other is a pre-image of $d_2$ or $d_3$, *semi-periodic*. The endpoint which is in $\overline{\mathcal{P}}_1$ will be referred to as an *aperiodic endpoint*.

Since there are finitely many pre-images of $d_2$ and $d_3$, the length $\Delta > 0$ of the shortest interval in either in the periodic set or the semi-periodic set is properly defined. We say that a point belongs to the *aperiodic set* if it does belong either to the periodic or the semi-periodic set.

**6. Behavior of the stochastic process when $\mathcal{P}$ is finite.** We can now show how the behavior of the stochastic process $\Xi$ introduced in Section 2 is influenced by the behavior of the deterministic triangle process. Recalling (9), we consider

$$\zeta(n) = \Lambda(\xi(\tau_n)),$$

the projection onto $\mathcal{S}_1^+$ of $\Xi$ at the random times $\tau_n, n = 0, 1, 2, \ldots$, when the server switches from one queue to another (where $\tau_0 = t_0$, a constant to be determined), under the transience conditions in (5), that is, $\lambda_i/\mu_i = \rho_i < 1$ for each $i$ and $\rho = \sum \rho_i > 1$.

We wish to prove Theorem 3.4, namely, that the trajectories of $\zeta(n)$ a.s. converge onto one of the stable orbits of the triangle process, under the assumption that $\mathcal{P}$, the set of pre-images of the decision points under $\varphi$, is finite and $\mathcal{P} \cap A^t = \varnothing$ for some finite $t$.

Recall the assumption that the service time distribution with the server at queue $i$ has finite variance $\sigma_i^2$. It is convenient, as with the triangle process, to use the transformation in Lemma 4.1, so we can assume that each $\mu_i = 1$. We now spell out some elementary properties of random walks and Poisson processes.

*Preliminaries.* Consider a random walk $S_n = \sum_{t=1}^n X_t$ on $\mathbb{Z}$, where $S_0 = 0$, $\mathbf{E}(X_t) = \rho_j - 1 < 0$ and $\operatorname{Var}(X_t) = \sigma^2$ and $X_t \geq -1$. Define the stopping times $T_c = \min\{n : S_n = -c\}$ for $c = 1, 2, \ldots$ and let $T_{1,n}$ be independent copies of $T_1$ so that $T_c$ is equal in distribution to $\sum_{t=1}^c T_{1,t}$. Using the standard martingales $S_n + n(\rho_j - 1)$ and $(S_n + n(\rho_j - 1))^2 - n\sigma^2$, we find that $\mathbf{E}(T_1) = 1/(1 - \rho_j)$ and $\operatorname{Var}(T_1) = \sigma^2/(1 - \rho_j)^3$ and, hence, that

$$\mathbf{E}(T_c) = c/(1 - \rho_j)$$

and

$$\operatorname{Var}(T_c) = c\sigma^2/(1 - \rho_j)^3.$$

Chebyshev's inequality provides the bound $\mathbf{P}[|T_c - c/(1 - \rho_j)| > c^{2/3}] \leq c\sigma^2/c^{4/3}(1 - \rho_j)^3 = c^{-1/3}\sigma^2/(1 - \rho_j)^3$ on the likely size of deviations from the mean.

Next consider a homogeneous Poisson process with arrivals at rate $\rho$ and let $L = L_0 + N(T_c)$, where $L_0$ is constant and $N(T_c)$ the number of arrivals by time $T_c$. Standard calculations give $\mathbf{E}(L) = L_0 + \rho\mathbf{E}(T_c)$ and $\mathrm{Var}(L) = \rho^2 \mathrm{Var}(T_c) + \rho\mathbf{E}(T_c)$. Applying these to the process $\Xi$ over the interval $\tau_{n+1} - \tau_n$ with the server at node $j$, we have, for $i \neq j$,

$$\mathbf{E}(\xi_i(\tau_{n+1})|\xi(\tau_n)) = \xi_i(\tau_n) + \frac{\rho_i \xi_j(\tau_n)}{(1-\rho_j)} \tag{17}$$

and

$$\mathrm{Var}(\xi_i(\tau_{n+1})|\xi(\tau_n)) = \xi_j(\tau_n) \frac{\rho_i^2 \sigma_j^2 + \rho_i(1-\rho_j)^2}{(1-\rho_j)^3}. \tag{18}$$

Using Chebyshev's inequality, we have

$$\begin{aligned} &\mathbf{P}\left[\left|\xi_i(\tau_{n+1}) - \xi_i(\tau_n) - \frac{\rho_i \xi_j(\tau_n)}{(1-\rho_j)}\right| > \xi_j(\tau_n)^{2/3} \middle| \xi(\tau_n)\right] \\ &\qquad \le \xi_j(\tau_n)^{-1/3} \frac{\rho_i^2 \sigma_j^2 + \rho_i(1-\rho_j)^2}{(1-\rho_j)^3}. \end{aligned} \tag{19}$$

While typical deviations on the number of arrivals are of order $\sqrt{\xi_j(\tau_n)}$, this simple bound for larger deviations will be enough for our purposes.

To estimate the effect of the projection, we apply Taylor's theorem to the function

$$g(u_1, u_2; a, b) = \frac{a + u_1}{b + u_1 + u_2},$$

where $0 \le a \le b$ are constants [compare with equation (10)], which provides the formula

$$g(u_1 + h_1, u_2 + h_2; a, b) - g(u_1, u_2; a, b) = \frac{h_1 u_2 - h_2(a + u_1)}{(b + u_1 + u_2 + \alpha(h_1 + h_2))^2},$$

where $\alpha \in (0, 1)$. Applying this at $u_1 + u_2 = u$ with $|h_i| \le u^{2/3}$ for $i = 1, 2$—the deviation considered in (19)—we obtain the bound

$$\begin{aligned} &|g(u_1 + h_1, u_2 + h_2; a, b) - g(u_1, u_2; a, b)| \\ &\qquad \le \frac{(a + u_1 + u_2)(|h_1| + |h_2|)}{(b + u_1 + u_2 + \alpha(h_1 + h_2))^2} \\ &\qquad \le \frac{|h_1| + |h_2|}{b + u_1 + u_2}\left(1 + 2\frac{|h_1| + |h_2|}{b + u_1 + u_2} + O(u^{-2/3})\right) \\ &\qquad \le 2u^{-1/3} \qquad \text{as } u \to \infty. \end{aligned} \tag{20}$$

Deviations $h_i = O(\sqrt{u}\,)$ give a smaller bound, but what we have will suffice.

PROOF OF THEOREM 3.4. We start with the stochastic process $W_n = \sum_{i=1}^3 \xi_i(\tau_n)$, $n = 0, 1, 2, \dots$, the total number of jobs at the queues at switching epochs subsequent to $\tau_0 = t_0$, some initial time to be chosen below. From (17) and (18), we have

$$\mathbf{E}(W_{n+1}|\xi(\tau_n), \mathfrak{R}(\xi(\tau_n)) = j) = \frac{\rho_i + \rho_k}{1-\rho_j}\xi_j(\tau_n) + \xi_i(\tau_n) + \xi_k(\tau_n)$$

$$= W_n + \frac{\theta}{1-\rho_j}\xi_j(\tau_n),$$

where $\theta = \sum_1^3 \rho_k - 1 > 0$. As we are only considering threshold switching rules, there is some value $D > 0$ (dependent on the decision points) such that if the server switches to queue $j$ at time $\tau_n$, then $\xi_j(\tau_n) > DW_n$ for each $j$ and every $n$. Hence, with $\nu = \frac{1}{2}\min_j D\theta/(1-\rho_j) > 0$, this implies that $\{(1+2\nu)^{-n}W_n\}$ is a nonnegative submartingale and so $W_n$ grows exponentially in mean. At this point we also introduce $\nu' > 0$ with $(1+\nu')^3 = 1+\nu$.

We are assuming that $\mathcal{P} \cap A^t = \varnothing$ eventually so, for any sufficiently small $\varepsilon_0 > 0$, there exists $n_0(\varepsilon_0)$ such that $|z - z'|_1 > 2\varepsilon_0$ whenever $z \in \mathcal{P}$, $z' \in A^{n_0}$, that is, the set $A^{n_0}$ is at least $2\varepsilon_0$ away from the decision points and their pre-images. Also, by Lemma 4.2, we can choose $n_0$ large enough that $A^{n_0} \subset C(\gamma)$, the set where $\varphi$ is contracting. We want to show that any stochastic trajectory (in projection) enters $A^{n_0}$ and then converges onto a trajectory of the triangle process.

Define a sequence of events $G_n$, $n = 0, 1, 2, \dots$, depending on parameters $w_0$, $\varepsilon_0$ by [recall (9)] $G_0(w_0) = \{W_0 \ge w_0\} \cap \{|\zeta(0) - z|_1 < \varepsilon_0 \text{ for some } z \in A^{n_0}\}$ and

$$G_n = \{W_n > w_0(1+\nu)^n\} \cap \{|\zeta(n) - \varphi(\zeta(n-1))|_1 \le \varepsilon_0\nu'(1+\nu')^{-n}\}.$$

We first bound $\mathbf{P}[G_n|\bigcap_1^{n-1} G_t]$. Choose $w_0$ large enough that $w^{1/3} > 1/\nu$ or, equivalently, $w^{2/3} < \nu w$ for all $w \ge w_0$. As $\mathbf{E}(W_n|W_{n-1}) \ge (1+2\nu)W_{n-1}$, we have from (19)

$$\begin{aligned}
&\mathbf{P}[W_n \le w_0(1+\nu)^n | W_{n-1} > w_0(1+\nu)^{n-1}] \\
&\quad < \mathbf{P}[|W_n - \mathbf{E}(W_n|W_{n-1})| > \nu w_0(1+\nu)^{n-1} | W_{n-1} > w_0(1+\nu)^{n-1}] \\
&\quad < \eta\nu^{-1/3}w_0^{-1/3}(1+\nu')^{-(n-1)},
\end{aligned}$$

where $\eta = \max_{ijk}((\rho_i+\rho_j)^2\sigma_k^2 + (\rho_i+\rho_j)(1-\rho_k)^2))/(1-\rho_k)^3$, a loose but adequate bound which shows that $W_n$ grows geometrically quickly with large probability.

Now we deal with deviation of the stochastic triangle process from the deterministic one. From any $y \in \mathbb{R}_+^3$ with $\Lambda(y) = z \in A^0 \setminus A_{\hat{\jmath}}^0$ with the server at node $\hat{\jmath}$, we have

$$f_{\hat{\jmath}}(z) = \sum_{i \ne \hat{\jmath}} \frac{(1-\rho_{\hat{\jmath}})z_i + \rho_i z_{\hat{\jmath}}}{(1-\rho_{\hat{\jmath}}) + \theta_{\hat{\jmath}} z_{\hat{\jmath}}} e_i \equiv \left(\frac{(1-\rho_{\hat{\jmath}})z_{\hat{\imath}} + \rho_{\hat{\imath}} z_{\hat{\jmath}}}{(1-\rho_{\hat{\jmath}})(z_{\hat{\imath}} + z_{\hat{k}}) + (\rho_{\hat{\imath}} + \rho_{\hat{k}})z_{\hat{\jmath}}}, \hat{\jmath}\right)$$

and multiplying the fraction by $\sum y_i/\sum y_i$ to switch to un-normalized values, we have

$$\frac{(1-\rho_{\hat{\jmath}})z_{\hat{\imath}}+\rho_{\hat{\imath}}z_{\hat{\jmath}}}{(1-\rho_{\hat{\jmath}})(z_{\hat{\imath}}+z_{\hat{k}})+(\rho_{\hat{\imath}}+\rho_{\hat{k}})z_{\hat{\jmath}}}=\frac{(1-\rho_{\hat{\jmath}})y_{\hat{\imath}}+\rho_{\hat{\imath}}y_{\hat{\jmath}}}{(1-\rho_{\hat{\jmath}})(y_{\hat{\imath}}+y_{\hat{k}})+(\rho_{\hat{\imath}}+\rho_{\hat{k}})y_{\hat{\jmath}}}. \tag{21}$$

For the stochastic process $\Xi$ starting from $\xi(\tau_{n-1})=y$ with $\sum_i y_i>w_0(1+\nu)^{n-1}$ and the server at $\hat{\jmath}$, the bound (19) implies that, for $i\neq\hat{\jmath}$,

$$\mathbf{P}\left[\left|\xi_i(\tau_n)-\left(y_i+\rho_i\frac{y_{\hat{\jmath}}}{1-\rho_{\hat{\jmath}}}\right)\right|>y_{\hat{\jmath}}^{2/3}\middle|\xi(\tau_{n-1})=y\right]\le\eta y_{\hat{\jmath}}^{-1/3}.$$

Now the identity (21) and the estimate (20), when applied at $\xi(\tau_{n-1})=y$ with $\sum y_i>w_0(1+\nu)^{n-1}$ and $\xi_i(\tau_n)$ in place of the $u_i+h_i$, immediately imply

$$\begin{aligned}\mathbf{P}[|\zeta(n)-\varphi(\zeta(n-1))|_1&<2w_0^{-1/3}(1+\nu')^{-(n-1)}|\xi(\tau_{n-1})]\\&\ge1-2\eta w_0^{-1/3}(1+\nu')^{-(n-1)}.\end{aligned}$$

This bound only relies upon $W_{n-1}=\sum_i\xi_i(\tau_{n-1})>w_0(1+\nu)^{n-1}$, an event which contains $\bigcap_1^{n-1}G_t$. As the process $\Xi$ is Markov and we can choose $w_0$ large enough that $2w_0^{-1/3}<\varepsilon_0\nu'/(1+\nu')$ and, as above, $w_0^{1/3}>1/\nu$, we have established, for such $w_0$, that

$$\mathbf{P}\left[G_n\middle|\bigcap_1^{n-1}G_t\right]\ge1-\eta\varepsilon_0\nu'(1+\nu')^{-n}.$$

A very similar argument can be used to show that $\mathbf{P}[G_0]\ge1-\eta\varepsilon_0$ when we choose $w_0$ such that $2n_0w_0^{-1/3}<\varepsilon_0$.

As the stochastic trajectory projects into $C(\gamma)$ (the region of $A^0$ where $\varphi$ is contracting), then, as long as the switching decision is the same at $\zeta(t)$ and $\varphi^{(t)}(\zeta(0))$ for every $t$, we have

$$\begin{aligned}\zeta(n)-\varphi^{(n)}(\zeta(0))&=\sum_{t=0}^{n-1}\varphi^{(t)}(\zeta(n-t))-\varphi^{(t+1)}(\zeta(n-t-1))\\&\le\sum_{t=0}^{n-1}\zeta(n-t)-\varphi(\zeta(n-t-1))\\&\le\varepsilon_0\nu'\sum_{t=1}^{n}(1+\nu')^{-t}<\varepsilon_0\qquad\text{on }\bigcap_0^nG_t.\end{aligned}$$

This inequality says that the trajectory $\zeta(t)$ is never further than $\varepsilon_0$ from the triangle process trajectory $z(t)$ started at $\zeta(0)$. By choice of $G_0$, every point of $z(t)$ is more than $\varepsilon_0$ from any $d_i$, so the switching decisions for $\zeta(t)$ and $z(t)$ will always be the same. This shows that on $\bigcap_0^\infty G_t$ the projections of the trajectories $\xi(t)$ converge a.s. onto deterministic trajectories of the triangle process which by Theorem 3.1 are periodic.

Finally,

$$\begin{aligned}\mathbf{P}\left[\bigcap_{t=0}^{n} G_t\right] &= \mathbf{P}[G_0]\prod_{t=1}^{n}\mathbf{P}\left[G_t\Big|\bigcap_{0}^{t-1} G_s\right]\\ &\ge \prod_{t=0}^{n-1}(1-\eta\varepsilon_0\nu'(1+\nu')^{-t})\\ &\to p(\varepsilon_0)>0 \qquad \text{as } n\to\infty,\end{aligned}$$

where $p(\varepsilon_0)\to 1$ as $\varepsilon_0\to 0$. This establishes Theorem 3.4, namely, that the stochastic trajectories $\zeta(t)$ converge onto the stable orbits of $z$, as we can choose $\varepsilon_0$ as small as we like. $\square$

Remark. The argument above is not delicate enough to decide whether trajectories $\zeta(t)$ converge with positive probability onto orbits which are stable on one side, for example, that are of even length and contain a decision point.

**7. The stochastic process when $\mathcal{P}$ is infinite.** The proofs of the following results employ the results of Section 5.3.

Lemma 7.1. *Let $\zeta(n)$ denote the state of the stochastic system at time $n$. Fix $\delta>0$. Define the event $E_{N,\delta}$ by*

$$E_{N,\delta}=\left\{\begin{array}{l}|\zeta(N+k)-\varphi^{(k)}(\zeta(N))|<\delta \text{ for all } k\ge 0 \text{ such that}\\ \mathrm{Side}(\zeta(N+i))=\mathrm{Side}(\varphi^{(i)}(\zeta(N))) \text{ for } i=0,1,\dots,k\end{array}\right\}. \tag{22}$$

*Then there exists a.s. $N=N(\delta,\omega)$ such that $E_{N,\delta}$ occurs.*

Proof. Analogously to Section 6, one can conclude that there is $0<\gamma<1$ such that a.s., for all large $n$,

$$\varepsilon_n := |\zeta(n)-\varphi(\zeta(n-1))|<\gamma^n. \tag{23}$$

Let $N_1=N_1(\omega)$ be the smallest of such $n$ and set $N:=\max(N_1,\min\{n: \sum_{i=n}^\infty \gamma^i<\delta\})$. Then as long as $\zeta(N+k)$ and $\varphi^{(k)}(\zeta(N))$ go through the same sequence of sides of the triangle $A^0$,

$$\begin{aligned}&|\zeta(N+k)-\varphi^{(k)}(\zeta(N))|\\ &\quad=|\zeta(N+k)-\varphi(\zeta(N+k-1))+\varphi(\zeta(N+k-1))-\varphi(\varphi^{(k-1)}(\zeta(N)))|\\ &\quad\le\gamma^{n+k}+|\zeta(N+k-1)-\varphi^{(k-1)}(\zeta(N))|,\end{aligned}$$

where we used the fact that $\varphi$ is contracting. Now the statement of the lemma follows from induction. □

Recall the definition of the aperiodic set and of $\Delta$ at the end of Section 5.3.

LEMMA 7.2. *Define the event $E$ by*

$$E = \{\textit{the stochastic system never comes closer than } \Delta/2 \textit{ to pre-images of } d_2 \textit{ or } d_3\}.$$

*Then, given that the stochastic system starts in the aperiodic set,* $\mathbf{P}(E) > 0$.

PROOF. First of all, we show that there is a constant $K > 0$ such that if $u$ belongs to the semi-periodic set, then $\varphi^{(k)}(u)$ belongs to the aperiodic set for some $k \le K$. Indeed, the image of an interval $I$ of the semi-periodic set must be either another semi-periodic interval, or must lie entirely inside an aperiodic interval or a semi-periodic interval. However, if $\varphi^{(k)}(I)$ is in one of the finitely many semi-periodic intervals for each $k$, this would imply that the points of $I$ are periodic. This contradicts Lemma 5.5.

Now let $\delta = \Delta/(2K)$ and suppose the system starts in the aperiodic set. By similar arguments to those of Section 6, with positive probability, $\sum_{n=0}^{\infty} |\varepsilon_n| < \delta$, where $\varepsilon_n$ is defined in (23). If for some $n$, $\zeta(n)$ is in the aperiodic set but $\zeta(n+1)$ is not, then $u := \zeta(n+1)$ cannot be "far" from the points of $\mathcal{P}_1$ since $|\varepsilon_n| < \delta$. Because the intervals of the periodic set are separated from the aperiodic set by the intervals of semi-periodic set, $u$ a.s. belongs to a semi-periodic interval, say, $I := [\tilde{l}, \tilde{r}]$, where $\tilde{l} \in \overline{\mathcal{P}}_1$ and $\tilde{r}$ is a pre-image of $d_2$ or $d_3$. We claim that, within the next $K$ applications of $\varphi$, $u$ will be "thrown out" of the image of $I$, or will be in the aperiodic set or some semi-periodic interval not far from its aperiodic endpoint again.

To show this, suppose $\zeta(n+k)$ lies in $\varphi^{(k)}(I)$ for the first $K$ steps. Then, by the choice of $K$, for some $k$, it will end up in the aperiodic set, and also it will not approach $\varphi^{(k)}(\tilde{r})$ (which *could* possibly be $d_2$ or $d_3$) closer than $K \times \delta \le \Delta/2$. On the other hand, if $\zeta(n+k-1) \in \varphi^{(k-1)}(I)$ and yet $\zeta(n+k) \notin \varphi^{(k)}(I)$, it means that $\zeta(n+k)$ is no further from $\varphi^{(k)}(\tilde{r})$ than $\varepsilon_{n+k}$ and it lies either in the aperiodic set or in a semi-periodic interval not far from its aperiodic endpoint $\varphi^{(k)}(\tilde{r})$. (Conditioned on $\sum_{n=0}^{\infty} |\varepsilon_n| < \delta$, there is not enough randomness for the stochastic system to exit the interval via the endpoint which is a pre-image of $d_2$ or $d_3$.) □

LEMMA 7.3. *Conditioned on the event $E$, periodicity of the stochastic system implies that $d_1$ is a limiting point for the trajectory of the stochastic system.*

Proof. Suppose not. Then there is $\delta = \delta(\omega) > 0$ such that the points in the $\delta$-neighborhood of $d_1$ are never hit by the stochastic system. Without loss of generality, suppose $\delta < \Delta/2$. By Lemma 7.1, there exists an $N$ for this $\delta$ such that (22) is fulfilled. Then for any $k \geq 0$, the distance between points $\varphi^{(k)}(\zeta(N))$ and $\zeta(N+k)$ does not exceed $\delta$ as long as they follow the same trajectory. On the other hand, they *will* follow the same trajectory, as there will never be the decision point $d_1$ between them, nor decision points $d_2$ or $d_3$, since we are conditioning on $E$.

Hence, $\tilde{l}(\varphi^{(k)}(\zeta(N)))$ will be also periodic, following the same path as $\varphi^{(k)}(\zeta(N))$, yielding contradiction with Lemma 5.3. □

Now we are able to finish the proof of Theorem 3.5.

Proof of Theorem 3.5. First, analogously to the proof of Theorem 3.4, one can easily show that, with positive probability, the stochastic system will converge to one of the finite cycles. So, it remains to show that, with positive probability, the stochastic system is not periodic.

Condition on the event $E$. Suppose that the stochastic system is periodic with period $m$, that is, for some $N_1$, we have

$$\mathrm{Side}(\zeta(n)) = \mathrm{Side}(\zeta(n+m)) \qquad \text{whenever } n \geq N_1. \tag{24}$$

For the moment consider only one image $\varphi(d_1) = f_2(d_1)$ of $d_1$. Since $f_2(d_1)$ is not periodic by Lemmas 5.4 and 5.3, there is a positive integer $K_1 = K_1(m) > m$ such that

$$\mathrm{Side}(\varphi^{(K_1)}(d_1)) \neq \mathrm{Side}(\varphi^{(K_1-m)}(d_1)). \tag{25}$$

Similarly, for the other possible image of $d_1$, $f_3(d_1)$, there is $K_2 = K_2(m) > m$ such that

$$\mathrm{Side}(\varphi^{(K_2)}(d_1)) \neq \mathrm{Side}(\varphi^{(K_2-m)}(d_1)), \tag{26}$$

under the assumption that $\varphi(d_1) \in A_3^0$.

Choose $\delta > 0$ so small, that each of the one-sided $\delta$-neighborhoods of $d_1$ which map onto sides $A_2^0$ and $A_3^0$, respectively, does not intersect with $d_1$ for the first $K = \max\{K_1, K_2\}$ applications of $\varphi$, and let $\delta_1 \leq \delta$ be the size of the smaller of these neighborhoods after $K$ mappings by $\varphi$.

By Lemma 7.3, the stochastic system will hit the $\delta_1/2$ neighborhood of $d_1$ at, say, time $N$. Since, in fact, the stochastic system will hit this neighborhood at arbitrary large times, we can suppose that $N > N_1$ and that $\sum_{n=N}^{\infty} |\varepsilon_n| < \delta_1/2$, where $\varepsilon_n$ is defined in (23). Also, for definiteness suppose that $\zeta(n)$ is on the side of $d_1$ which maps onto $A_2^0$. Then the stochastic system will follow the image of $d_1$ which maps onto $A_2^0$ for the next $K$ steps, in the sense $\mathrm{Side}(\varphi^{(k)}(\zeta(n))) = \mathrm{Side}(\varphi^{(k)}(d_1))$ since $\delta_1/2 + \delta_1/2 \leq \delta_1$. However, recall that $K$ was chosen in such a way that the sequence $\mathrm{Side}(\varphi^{(k)}(d_1))$, $k = 1, 2, \ldots, K$, cannot be $m$-periodic, creating the contradiction between (24) and (25). □

## REFERENCES

[1] Asmussen, S. (2003). *Applied Probability and Queues*, 2nd ed. Springer, New York. MR1978607

[2] Coelho, Z., Lopes, A. and da Rocha, L. F. (1995). Absolutely continuous invariant measures for a class of affine interval exchange maps. *Proc. Amer. Math. Soc.* **123** 3533–3542. MR1322918

[3] Dajani, K. (2002). A note on rotations and interval exchange transformations on 3-intervals. *Int. J. Pure Appl. Math.***1** 151–160. MR1900223

[4] Fayolle, G., Malyshev, V. A. and Men′shikov, M. V. (1995). *Topics in the Constructive Theory of Countable Markov Chains*. Cambridge Univ. Press. MR1331145

[5] Foss, S. and Last, G. (1996). Stability of polling systems with exhaustive service policies and state-dependent routing. *Ann. Appl. Probab.* **6** 116–137. MR1389834

[6] Foss, S. and Last, G. (1998). On the stability of greedy polling systems with general service policies. *Probab. Engrg. Inform. Sci.* **12** 49–68. MR1492140

[7] Liousse, I. and Marzougui, H. (2002). Échanges d'intervalles affines conjugués à des linéaires. *Ergodic Theory Dynam. Systems* **22** 535–554. MR1898804

[8] MacPhee, I. M. and Menshikov, M. V. (2003). Critical random walks on two-dimensional complexes with applications to polling systems. *Ann. Appl. Probab.* **13** 1399–1422. MR2023881

[9] Menshikov, M. and Zuyev, S. (2001). Polling systems in the critical regime. *Stochastic Process. Appl.* **92** 201–218. MR1817586

[10] Meyn, S. P. and Tweedie, R. L. (1993). *Markov Chains and Stochastic Stability*. Springer, London. MR1287609

[11] Veech, W. A. (1984). The metric theory of interval exchange transformations, I–III. *Amer. J. Math.* **106** 1331–1422. MR0765582

I. M. MacPhee
M. V. Menshikov
Department of Mathematical Sciences
University of Durham Science Site
Durham DH1 3LE
United Kingdom
E-mail: i.m.macphee@durham.ac.uk
mikhail.menshikov@durham.ac.uk

S. Popov
Instituto de Matemática e Estatística
Universidade de São Paulo
São Paulo SP
Brasil
E-mail: popov@ime.usp.br

S. Volkov
School of Mathematics
University of Bristol
Bristol BS8 1TW
United Kingdom
E-mail: s.volkov@bristol.ac.uk